\documentclass[12pt,a4paper]{smfart}
\usepackage{amsmath,amsthm,upref}

\numberwithin{equation}{section}

\usepackage[dvips]{graphicx}
\usepackage{psfrag}
\usepackage{a4wide}
\usepackage[francais]{babel}
\usepackage{epsfig,subfigure}

\newcommand\Su{\mathcal S}
\newcommand\Sl{\textrm{SL}_2(\mathbb R)}
\newcommand\C{\mathbb C}
\newcommand\R{\mathbb R}
\newcommand\N{\mathbb{N}}
\newcommand\Z{\mathbb Z}
\newcommand\T{\mathcal T}
\newcommand\F{\mathcal F}

\newcommand\D{{\rm d}\hspace{-0.5pt}}
\DeclareMathOperator{\Leb}{Leb}
\DeclareMathOperator{\Aire}{Aire}
\DeclareMathOperator{\sys}{sys}
\DeclareMathOperator{\Card}{Card}
\DeclareMathOperator{\Comp}{Comp}
\DeclareMathOperator{\Int}{Int}

\newcommand{\Sbb}{\mathbb{S}}
\newcommand{\numero}{\subsection{}}
\newcommand{\tq}{\; : \; }
\newcommand{\moins}{-}
\newcommand{\A}{\mathcal{A}}
\newcommand{\M}{\mathcal{M}}
\renewcommand{\epsilon}{\varepsilon}

\swapnumbers
\newtheorem{Theorem}[subsection]{Th\'{e}or\`{e}me}
\newtheorem{Corollary}[subsection]{Corollaire}
\newtheorem{Proposition}[subsection]{Proposition}
\newtheorem{Lemma}[subsection]{Lemme}
\newtheorem*{NoNumberTheorem}{Th\'{e}or\`{e}me}
\newtheorem*{NoNumberProposition}{Proposition}

\newtheorem*{ThmA}{Th\'{e}or\`{e}me A}
\newtheorem*{ThmB}{Th\'{e}or\`{e}me B}
\newtheorem{Definition}[subsection]{D\'{e}finition}
\theoremstyle{remark}
\newtheorem{Remark}[subsection]{Remarque}

\begin{document}
\title[Unique ergodicit\'{e} sur les surfaces plates]
{Un th\'{e}or\`{e}me de Kerckhoff, Masur et Smillie~: \\
Unique ergodicit\'{e} sur les surfaces plates}

\author{S\'{e}bastien Gou\"{e}zel}

\address{
IRMAR, Universit\'{e} de Rennes 1, Campus de Beaulieu, B\^{a}timent 22,
35042 Rennes Cedex, France. }

\email{sebastien.gouezel@univ-rennes1.fr}

\urladdr{perso.univ-rennes1.fr/sebastien.gouezel}

\author{Erwan Lanneau}

\address{
Centre de Physique Th\'{e}orique de Marseille (CPT), UMR CNRS 6207
\newline Universit\'{e} du Sud Toulon-Var et \newline F\'{e}d\'{e}ration de
Recherches des Unit\'{e}s de Math\'{e}matiques de Marseille \newline Luminy,
Case 907, F-13288 Marseille Cedex 9, France. }

\email{lanneau@cpt.univ-mrs.fr}

\urladdr{www.cpt.univ-mrs.fr/~lanneau}


\date{27 novembre 2006}

\begin{abstract}
Ces notes correspondent \`{a} un cours donn\'{e} lors de l'\'{e}cole th\'{e}matique
de th\'{e}orie ergodique au C.I.R.M. \`{a} Marseille
en avril 2006. Nous pr\'{e}sentons et d\'{e}montrons un
th\'{e}or\`{e}me de Kerckhoff, Masur et Smillie sur l'unique ergodicit\'{e} du
flot directionnel sur une surface de translation dans presque toutes
les directions. \\
La preuve suit essentiellement celle pr\'{e}sent\'{e}e dans un survol
de Masur et Tabachnikov. Nous donnons une preuve compl\`{e}te et
\'{e}l\'{e}mentaire du th\'{e}or\`{e}me.
\end{abstract}

\maketitle

\setcounter{tocdepth}{1}

\tableofcontents
\newpage
\section{Introduction}

Le but de ce texte est de d\'{e}montrer un th\'{e}or\`{e}me de Steven
Kerckhoff, Howard Masur et John Smillie concernant un r\'{e}sultat de
th\'{e}orie ergodique sur les surfaces de translation
(\cite{Kerckhoff:Masur:Smillie}).

\begin{figure}[htbp]
\begin{center}

\psfrag{f}{$\scriptstyle \F_\theta^t(x)$} \psfrag{x}{$\scriptstyle
x$}  \psfrag{T}{$\scriptstyle \mathbb T^2$}

  \subfigure[]{\epsfig{figure=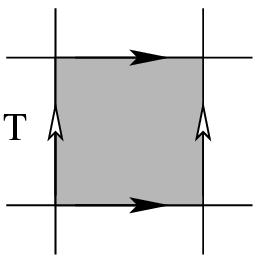,width=3.0cm}} \qquad \qquad
  \subfigure[]{\epsfig{figure=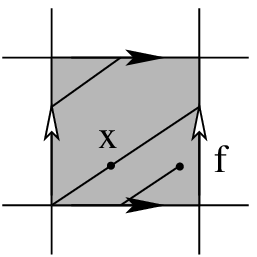,width=3.0cm}}

\end{center}
\caption{\label{fig:tore} Le tore $\R^2 / \Z^2$ muni de son flot
directionnel $\F_\theta$. }
\end{figure}

\numero Le prototype d'une surface de translation est donn\'{e} par
le tore plat standard $\mathbb T^2=\R^2 /\Z^2$ muni de sa
m\'{e}trique plate (ou euclidienne) $\D z$. On peut construire
g\'{e}om\'{e}triquement $\mathbb T^2$ de la mani\`{e}re suivante. Un
domaine fondamental pour l'action de $\mathbb Z^2$ sur $\mathbb
R^2$ par translation est donn\'{e} par le carr\'{e} unit\'{e} $]0,1[^2$. Le
tore s'obtient alors en identifiant les bords oppos\'{e}s de ce
carr\'{e} \`{a} l'aide de translations (voir figure~\ref{fig:tore}a).

\numero On peut d\'{e}finir sur $\mathbb T^2$ un flot directionnel
$\mathcal F_\theta$, pour $\theta \in \mathbb S^1$. C'est le
flot lin\'{e}aire de pente constante \'{e}gale \`{a} $\theta$~: si $x\in
\mathbb T^2$ et $t>0$ alors $\mathcal F_\theta^t(x)$ est
l'unique point de $\mathbb T^2$ \`{a} distance $t$ de $x$ et situ\'{e}
sur la g\'{e}od\'{e}sique (orient\'{e}e positivement) de pente $\theta$
passant par $x$ (voir figure~\ref{fig:tore}b). On d\'{e}finit de
mani\`{e}re \'{e}quivalente $\F_\theta^t(x)$ pour $t<0$.

Par construction il est clair que le flot $\mathcal F_\theta$ laisse
invariant la mesure de Lebesgue not\'{e}e $\Leb$. On peut donc
consid\'{e}rer le syst\`{e}me dynamique $(\mathbb T^2,\Leb,\mathcal
F_\theta)$ ; il v\'{e}rifie le th\'{e}or\`{e}me classique suivant,
d\^{u} \`{a} Hermann Weyl~:

\begin{NoNumberTheorem}
Pour Lebesgue presque tout $\theta\in\mathbb S^1$, le syst\`{e}me
dynamique $(\mathbb T^2,\Leb,\mathcal F_\theta)$ est uniquement
ergodique.
\end{NoNumberTheorem}
Rappelons que cela signifie que $\mathcal F_\theta$ admet une
seule mesure de probabilit\'{e} invariante. De mani\`{e}re \'{e}quivalente,
toutes les orbites du flot sont denses dans le tore et
uniform\'{e}ment distribu\'{e}es par rapport \`{a} la mesure de Lebesgue.
On dira aussi souvent dans la suite que le flot directionnel
est uniquement ergodique dans presque toutes les directions. En
fait, le th\'{e}or\`{e}me de Weyl est m\^{e}me plus pr\'{e}cis que cela,
puisqu'il donne l'unique ergodicit\'{e} dans toutes les directions
\emph{irrationnelles}~; les directions exceptionnelles forment
donc un ensemble d\'{e}nombrable.

\numero Une surface de translation est un triplet
$(\Su,\Sigma,\omega)$ o\`{u} $\Su$ est une surface compacte, connexe,
sans bord, orient\'{e}e, $\Sigma=\{P_1,\dots,P_n\}$ est un ensemble fini
de points $P_i\in \Su$ et $\omega=\left\{ (U_i,z_i) \right\}_i$ est
un atlas de translation sur $\Su\setminus \Sigma$. Par atlas de
translation, nous demandons que les changements de cartes soient du type
$z_i=z_j+cst$. On demandera de plus que, au voisinage de chaque
singularit\'{e} $P_i$, $\Sigma$ soit isomorphe \`{a} un rev\^{e}tement de
$\R^2 \moins \{0\}$ avec un nombre fini $k_i+1$ de feuillets.
Par abus de langage nous
d\'{e}noterons souvent $(\Su,\Sigma,\omega)$
simplement par $(\Su,\omega)$ ou juste~$\Su$ lorsque le contexte sera
clair. La structure euclidienne sur $\Su$ induit naturellement une
mesure que l'on appellera encore mesure de Lebesgue $\Leb$.

\numero Comme pr\'{e}c\'{e}demment, nous pouvons d\'{e}finir un flot
directionnel $\mathcal F_\theta$ sur $\Su$ laissant invariant la
m\'{e}trique euclidienne d\'{e}finie sur $\Su$, via les cartes. Ce flot
n'est en fait pas correctement d\'{e}fini sur $\Sigma$, ni m\^{e}me sur les
trajectoires qui arrivent ou partent de $\Sigma$ (il y en a un
nombre fini). Ainsi, $\mathcal F_\theta$ est d\'{e}fini sur un
sous-ensemble dense et de mesure pleine de $\Su$ mais, par abus de
langage, nous parlerons n\'{e}anmoins du flot $\F_\theta$ sur $\Su$.
Comme le flot directionnel dans $\R^2$ pr\'{e}serve la mesure de
Lebesgue, $\F_\theta$ pr\'{e}serve la mesure $\Leb$ sur $\Su$.

On peut alors \'{e}noncer le th\'{e}or\`{e}me que nous allons
d\'{e}montrer dans ce texte et qui g\'{e}n\'{e}ralise le
th\'{e}or\`{e}me de Weyl \'{e}nonc\'{e} ci-dessus.

\begin{NoNumberTheorem}[S.~Kerckhoff, H.~Masur, J.~Smillie, 1986]
Pour toute surface de translation $\Su$ et pour Lebesgue presque
tout $\theta\in\mathbb S^1$, le flot directionnel $\mathcal
F_\theta$ sur $\Su$ est uniquement ergodique.
\end{NoNumberTheorem}

\begin{Remark}
Le th\'{e}or\`{e}me ci-dessus reste vrai dans un cadre plus
g\'{e}n\'{e}ral. Nous pouvons \'{e}tendre la d\'{e}finition de
surfaces de translation aux surfaces de demi-translation en imposant
que les changements de cartes soient de la forme $z_i=\pm z_j +
cst$. L'objet naturel venant avec un atlas de translation est une forme
diff\'{e}rentielle ab\'{e}lienne. L'analogue pour les surfaces de demi-translation
sont les formes diff\'{e}rentielles quadratiques
(voir~\cite{Hubbard:Masur} pour un expos\'{e} plus d\'{e}taill\'{e}).

La preuve du th\'{e}or\`{e}me de Kerckhoff, Masur et Smillie dans ce cadre
\'{e}tendu est simplement plus technique et ne n\'{e}cessite pas d'id\'{e}e
nouvelle fondamentale. Par cons\'{e}quent, nous nous restreindrons aux
surfaces de translation.
\end{Remark}

\begin{Remark}
Remarquons aussi que, contrairement au tore, il y a des
exemples de surfaces de translation o\`{u} il y a un nombre non
d\'{e}nombrable de directions minimales et non uniquement
ergodiques (\cite{Veech:68,Keane}). Une des diff\'{e}rences
essentielles avec le tore est que le flot lin\'{e}aire sur un tore
est une isom\'{e}trie globale, i.e. la m\'etrique plate ne poss\`ede pas 
de singularit\'es.
\end{Remark}

\begin{Remark}
Enfin on peut noter que le th\'{e}or\`{e}me KMS est un renforcement
d'un th\'{e}or\`{e}me de H.~Masur et W.~Veech
(\cite{Masur:82,Veech:82}). Une section de Poincar\'{e} du flot
directionnel sur un intervalle dans $\Su$ produit un {\it
\'{e}change d'intervalles}. Masur et Veech ont d\'{e}montr\'{e} que presque
tout \'{e}change d'intervalles est uniquement ergodique.
\end{Remark}

\numero Dans cette note, nous allons suivre une d\'{e}marche un peu
atypique~:
apr\`{e}s avoir donn\'{e} quelques exemples et propri\'{e}t\'{e}s basiques des
surfaces de translation, et d\'{e}montr\'{e} une version faible du
th\'{e}or\`{e}me KMS, nous exposerons une preuve compl\`{e}te et
\'{e}l\'{e}mentaire du th\'{e}or\`{e}me KMS. \\
Ce n'est qu'ensuite que nous
introduirons des outils conceptuels suppl\'{e}mentaires qui apporteront
un autre \'{e}clairage sur la preuve. Nous esp\'{e}rons que la
d\'{e}monstration du th\'{e}or\`{e}me pourra \^{e}tre vue comme une
justification \emph{a priori} de l'introduction de ces outils.

\section{Quelques d\'{e}finitions et exemples}

\subsection{Surfaces de translation}

Avant d'aller plus loin, donnons quelques propri\'{e}t\'{e}s utiles des
surfaces de translation. Comme nous l'avons vu pr\'{e}c\'{e}demment,
$\Su\setminus \Sigma$ poss\`{e}de un atlas de translation. Ceci permet
de d\'{e}finir une m\'{e}trique euclidienne via la forme diff\'{e}rentielle
globale $\omega=\D z_i$ sur $U_i$. Ainsi $\Su\setminus \Sigma$ est
localement isom\'{e}trique \`{a} un plan $\mathbb R^2$, la courbure de la
m\'{e}trique \'{e}tant nulle en tout point. La formule de Gau{\ss}-Bonnet
implique alors que la courbure de la m\'{e}trique sur $\Su$ est
concentr\'{e}e dans les points $P_i$. \medskip

\noindent Par d\'{e}finition, il existe pour tout $i$ une identification
entre un voisinage \'{e}point\'{e} de $P_i$ dans $\Su$ et un
rev\^{e}tement de degr\'{e} $k_i+1$ de $\R^2 \moins \{0\}$.
La surface $\Su$ est
ainsi isom\'{e}trique \`{a} un c\^{o}ne d'angle $2 (k_i+1)\pi$. Nous
parlerons ainsi pour $P_i$ de singularit\'{e} conique d'angle $2 (k_i+1)\pi$.
Notons ici qu'au voisinage d'un point r\'{e}gulier, la surface est bien
isom\'{e}trique \`{a} un ``c\^{o}ne plat'' d'angle $2 (k+1)\pi = 2\pi$ avec
$k=0$. La courbure de Gau{\ss} est donn\'{e}e par $\kappa=-k_i\pi$ et la
formule de Gau{\ss}-Bonnet s'\'{e}crit ici~:
  \begin{equation}
  \sum_{i=1}^n k_i = 2g-2 \qquad \textrm{o\`{u} } g=\textrm{genre}(\Su).
  \end{equation}
Le th\'{e}or\`{e}me de Riemann-Roch implique que pour toute partition
enti\`{e}re $(k_1,\dots,k_n)$ de $2g-2$, il existe une surface de
translation avec exactement $n$ singularit\'{e}s coniques d'angles $2
(k_i+1)\pi$ pour $i=1,\dots,n$~; mais nous n'utiliserons pas ce fait
ici. \medskip

\noindent On peut d\'{e}finir un feuilletage vertical (respectivement
horizontal) sur $\Su\setminus \Sigma$ par les lignes de niveaux des cartes
locales~: $z_i^{-1}(x=cst)$ (respectivement $z_i^{-1}(y=cst)$). Cela
permet alors, dans ces coordonn\'{e}es, d'introduire le flot
directionnel $\mathcal F_\theta$ pour $\theta\in \mathbb S^1$. La
m\'{e}trique euclidienne est invariante sous l'action de ce flot. Les
orbites du flot directionnel ne sont bien s\^{u}r pas toutes bien
d\'{e}finies~: certaines feuilles rencontrent des singularit\'{e}s en
temps fini. N\'{e}anmoins $\Su$ ne poss\`{e}de qu'un nombre fini de
singularit\'{e}s et par chaque singularit\'{e} ne passe qu'un nombre
fini de feuilles dans la direction $\theta$. Ceci implique alors que
pour tout $\theta$, le flot $\F_\theta$ est presque partout bien
d\'{e}fini sur $\Su$. \medskip

\noindent La structure euclidienne sur $\Su$ est tr\`{e}s facile \`{a}
visualiser localement~; dans les coordonn\'{e}es fournies par le
feuilletage horizontal et le feuilletage vertical, une
g\'{e}od\'{e}sique est localement une droite de pente constante. Bien
s\^{u}r globalement une feuille peut avoir des comportements
baroques. Dans toute la suite, sauf mention du contraire, nous
allons nous restreindre aux surfaces d'aire $1$, c'est-\`{a}-dire
$\Aire(S) := \int_S \omega \wedge \overline{\omega} = 1$.
\medskip

\noindent Nous utiliserons les terminologies suivantes pour le flot
directionnel~:
\begin{itemize}
\item Une feuille passant par une singularit\'{e} est une \emph{s\'{e}paratrice}.
\item Une feuille connectant deux singularit\'{e}s (\'{e}ventuellement les
m\^{e}mes) est une \emph{connexion de selles} ou segment g\'{e}od\'{e}sique.
\item Une feuille ne passant pas par une singularit\'{e} est dite
  \emph{r\'{e}guli\`{e}re}.
\end{itemize}

\numero
\label{par:RecollePolygones}
La g\'{e}om\'{e}trie des surfaces de translation peut se comprendre en
voyant une telle surface comme recollement de polygones. Plus
pr\'{e}cis\'{e}ment, consid\'{e}rons des polygones $P_1,\dots, P_n$ du plan,
d'adh\'{e}rences disjointes, tels que les c\^{o}t\'{e}s de ces polygones
soient regroup\'{e}s par paires de c\^{o}t\'{e}s parall\`{e}les et de m\^{e}me
longueur. On peut alors former une surface $\Su$ en identifiant les c\^{o}t\'{e}s
correspondants. Nous ajoutons la conditions que l'angle conique autour des
sommets des polygones (vu dans $\Su$) soit de la forme $2c\pi$ o\`{u} $c$ est
un entier (sans cette condition, la notion d'angle est mal d\'efinie, 
voir par exemple la d\'efinition analogue dans~\cite{Masur:06}). 
La surface $\Su$ devient alors une surface de
translation.
C'est une surface de translation, dont les singularit\'{e}s forment
un sous-ensemble des sommets des polygones $P_j$.
\medskip

\noindent R\'{e}ciproquement, toute surface de translation peut
s'obtenir ainsi, comme nous allons le voir. Notons aussi que pour des
raisons \'{e}videntes, nous ne consid\'{e}rerons que des surfaces connexes.

\begin{Definition}
\label{def:Triangulation}
On appelle triangulation d'une surface de translation une
triangulation dont les sommets sont les singularit\'{e}s et les ar\^{e}tes
des connexions de selles.
\end{Definition}

\numero
\label{ExisteTriangulation}
Si on part d'un ensemble de connexions de selles dont les int\'{e}rieurs
sont deux \`{a} deux disjoints, on peut les compl\'{e}ter pour obtenir une
triangulation de la surface (ce qui se v\'{e}rifie ais\'{e}ment en
ajoutant des ar\^{e}tes tant que c'est possible et en remarquant que
quand ce n'est plus possible la surface est n\'{e}cessairement triangul\'{e}e).

\noindent Ainsi, toute surface admet une triangulation. On en d\'{e}duit
ais\'{e}ment que toute surface de translation s'obtient en recollant des
triangles suivant la proc\'{e}dure d\'{e}crite en \ref{par:RecollePolygones}.

\numero
Consid\'{e}rons une triangulation d'une surface de translation $\Su$ de
genre $g$, et notons respectivement $s,e,f$ le nombre de sommets,
c\^{o}t\'{e}s et faces de cette triangulation. Alors $s-e+f=2-2g$ par la
formule d'Euler. Comme chaque face a trois c\^{o}t\'{e}s et chaque
c\^{o}t\'{e} est sur le bord de deux faces, n\'{e}cessairement
$2e=3f$. Comme $s$ est le nombre de singularit\'{e}s, on obtient
  \begin{equation}
  \label{FormuleFaces}
  f=2(s+2g-2).
  \end{equation}

\subsection{Exemples}

\begin{figure}[htbp]
\begin{center}
\psfrag{a}{$6\pi$} \psfrag{S}{$\Su$}

\includegraphics[width=9cm]{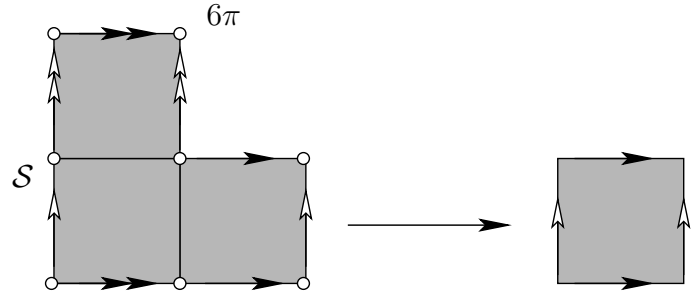}
\end{center}
\caption{ \label{fig:exemple:3:carreaux} Un rev\^{e}tement ramifi\'{e} de
degr\'{e} $3$ au dessus du tore. Le rev\^{e}tement est ramifi\'{e} au dessus de
l'origine et l'unique point critique est une singularit\'{e} pour la
m\'{e}trique, d'angle $6\pi$. On a donc $2\cdot \textrm{genre}(\Su)-2=2$, soit
$\textrm{genre}(\Su)=2$.
}
\end{figure}

Le tore plat fournit beaucoup d'exemples d\'{e}riv\'{e}s. Soit $f : \Su
\rightarrow \mathbb T^2$ un rev\^{e}tement ramifi\'{e} au dessus de
l'origine. Si on note $\D z$ la m\'{e}trique plate sur $\mathbb T^2$
alors la m\'{e}trique $f^\ast\D z$ \'{e}quipe $\Su\setminus f^{-1}(0)$ d'une
structure de surface plate (voir figure~\ref{fig:exemple:3:carreaux}).

\begin{figure}[htbp]
\begin{center}
\psfrag{a}{$6\pi$} \psfrag{S}{$\Su$}

\includegraphics[width=4.5cm]{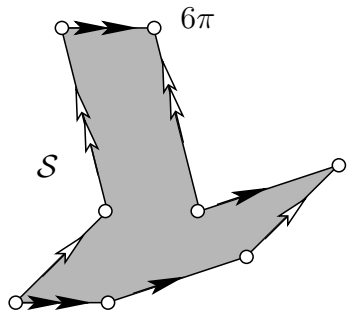}
\end{center}
\caption{ \label{fig:exemple:3:carreaux:modif} Une  surface de
translation de genre $2$. L'angle conique de l'unique singularit\'{e}
(point blanc) est $6\pi$.}
\end{figure}

\noindent Un exemple beaucoup plus ``g\'{e}n\'{e}rique'' peut \^{e}tre
construit \`{a} partir d'une surface plate quelconque, obtenue en
recollant des polygones~: il suffit de bouger les c\^{o}t\'{e}s des
polygones de mani\`{e}re \`{a} pr\'{e}server la propri\'{e}t\'{e} des paires de
c\^{o}t\'{e}s (la figure~\ref{fig:exemple:3:carreaux:modif} est un
exemple ``d\'{e}form\'{e}'' de la figure~\ref{fig:exemple:3:carreaux}).
Le terme g\'{e}n\'{e}rique employ\'{e} ci-dessus sera expliqu\'{e} dans la
section~\ref{sec:espace:module}.

\begin{figure}[htbp]
\begin{center}
\psfrag{a}{$\scriptstyle a$}  \psfrag{b}{$\scriptstyle b$}
\psfrag{c}{$\scriptstyle c$}  \psfrag{d}{$\scriptstyle d$}
\psfrag{e}{$\scriptstyle e$} \psfrag{cp}{$\scriptstyle c'$}
\psfrag{dp}{$\scriptstyle d'$} \psfrag{ep}{$\scriptstyle e'$}
 \psfrag{S}{$\Su_\lambda$}
\psfrag{P1}{$\scriptstyle P_1$}  \psfrag{P2}{$\scriptstyle P_2$}

\includegraphics[width=6.5cm]{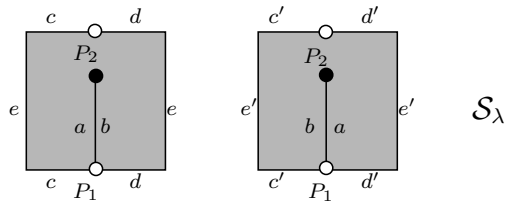}
\end{center}
\caption{ \label{fig:exemple:4:slit} Deux tores unit\'{e} recoll\'{e}s
le long du segment g\'{e}od\'{e}sique vertical $P_1P_2$ centr\'{e} et de
hauteur $\lambda \in ]0,1[$ (on identifie les c\^{o}t\'{e}s de m\^{e}me
\'{e}tiquette par translation). La surface $\Su_\lambda$
correspondante est de genre $2$. Elle poss\`{e}de deux singularit\'{e}s
coniques (points blanc et noir) chacune d'angle $4\pi$. }
\end{figure}

\noindent Nous donnons dans la figure~\ref{fig:exemple:4:slit}
un dernier exemple, toujours de genre $2$, mais cette fois-ci
avec deux singularit\'{e}s coniques, n\'{e}cessairement chacune d'angle
$4\pi$.

\subsection{Billards rationnels et surfaces de translation}

Soit $P \subset \R^2$ un polygone. Dans tout ce paragraphe, on
notera $\Gamma$ le sous groupe des isom\'{e}tries lin\'{e}aires de
$\R^2$ engendr\'{e} par les parties lin\'{e}aires des r\'{e}flexions par
rapport aux c\^{o}t\'{e}s de $P$.

\noindent On note $\tilde \F$ le flot du billard sur $P\times \Sbb^1$~: le
point $\tilde \F_t(x,\theta)$ est obtenu en partant de $x$, en
suivant une ligne droite dans la direction $\theta$ jusqu'\`{a}
rencontrer un c\^{o}t\'{e} du polygone $P$. On ``rebondit'' alors suivant
les lois de Descartes de l'optique g\'{e}om\'{e}trique
et on continue ce processus pendant un temps $t$.
Ce flot pr\'{e}serve la mesure $\Leb_P \otimes
\Leb_{\Sbb^1}$, et il est d\'{e}fini pour toutes les trajectoires qui
n'arrivent pas dans un des coins du polygone. Lorsque les angles du
billard sont des multiples rationnels de $\pi$, ce flot n'est
\'{e}videmment pas ergodique~: si $\theta \in \Sbb^1$, l'ensemble
$P\times \bigcup_{\gamma\in \Gamma} \{\gamma \cdot \theta\}$ est
invariant par $\tilde \F$ (notons que $\Gamma$ est alors fini). On
notera alors $\tilde \F_\theta$ la restriction de $\tilde \F$ \`{a}
cet ensemble~; ce flot pr\'{e}serve la mesure $\Leb \otimes
\sum_{\gamma \in \Gamma} \delta_{\gamma\cdot\theta}$. \medskip

\noindent Une construction classique (voir~\cite{FoKe,KaZa}) permet d'associer
\`{a} chaque polygone $P$ une surface de translation (non n\'{e}cessairement
compacte) $\Su(P)$ telle que le flot du billard et le flot
directionnel commutent avec la construction. Nous r\'{e}sumons ici la
construction en question. \medskip

\noindent Pour $\gamma \in \Gamma$ on notera
$P_\gamma=\gamma(P)$. Nous allons d\'{e}crire une relation
d'\'{e}quivalence sur la r\'{e}union disjointe $\bigsqcup_{\gamma \in
\Gamma} P_\gamma$ des polygones $P_\gamma$. Si $c$ est un c\^{o}t\'{e}
de $P_\gamma$ soit $\delta_c\in\Gamma$ la partie lin\'{e}aire de la
r\'{e}flexion  par rapport \`{a} $c$. On identifie alors $c\subset
P_\gamma$ avec $\delta_c(c)\subset P_{\delta \circ \gamma}$ par
translation. L'ensemble quotient pour cette relation est not\'{e}
$\Su(P)$. \medskip

\noindent La surface $\Su(P) \setminus \{ \textrm{sommets de } P \}$
h\'{e}rite naturellement d'un atlas de translation via les polygones
$P_\gamma$. Il est facile de v\'{e}rifier que cet atlas s'\'{e}tend
globalement en une structure de translation sur $\Su$ tout entier. Les
singularit\'{e}s de la m\'{e}trique forment un sous-ensemble
(\'{e}ventuellement strict) des sommets des
polygones $P_\gamma$. Par ailleurs, la surface $\Su$ est compacte si
et seulement si le polygone $P$ est rationnel, c'est \`{a} dire si et
seulement si les angles de $P$ sont des multiples rationnels de
$\pi$. Le polygone $P$ et ses images par les \'{e}l\'{e}ments de $\Gamma$
pavent alors la surface $\Su$ et on obtient une application $\pi :
\Su(P) \rightarrow P$. Par construction les orbites du flot $\F$ se
projettent sur les orbites du flot $\tilde \F$. Lorsque $\Su$ est
compacte, les orbites de $\F_\theta$ se projettent m\^{e}me sur celles
de $\tilde\F_\theta$, de mani\`{e}re $\pi$-\'{e}quivariante.

\noindent Ainsi comme corollaire direct du th\'{e}or\`{e}me KMS on
obtient le r\'{e}sultat suivant.

\begin{Corollary}
\label{cor:KMS_billards}
Pour tout billard rationnel, l'ensemble
  \begin{equation*}
  \{\theta \in \mathbb S^1,\ \tilde \F_\theta \textrm{ est uniquement
  ergodique} \}
  \end{equation*}
est de mesure pleine.
\end{Corollary}

\noindent Notons que ce r\'{e}sultat a des cons\'{e}quences pour les
billards non n\'{e}cessairement rationnels. En effet, en approchant un
billard (polygonal) arbitraire par des billards rationnels, gr\^{a}ce \`{a} une
remarque de A.~Katok, il permet de d\'{e}montrer le r\'{e}sultat suivant.

\begin{Corollary}
\label{cor:ergodique}
L'ensemble $\{ Q\ ; \ Q \textrm{ est un polygone \`{a} } k \textrm{
  c\^{o}t\'{e}s et } \tilde \F \textrm{ est ergodique} \}$ forme un
$G_\delta$ dense dans l'ensemble des polygones \`{a} $k$ c\^{o}t\'{e}s
(pour la topologie produit, c'est \`{a} dire en identifiant
l'ensemble des polygones \`{a} $k$ c\^{o}t\'{e}s avec un ouvert de
$\R^{2k}$).
\end{Corollary}

\numero Si l'on est principalement int\'{e}ress\'{e} par les billards, et en
particulier par le corollaire~\ref{cor:KMS_billards}, on pourrait
essayer de ne pas mentionner du tout les surfaces de translation, et
de transcrire directement la preuve du th\'{e}or\`{e}me KMS que nous
allons pr\'{e}senter dans le langage des billards. Cependant, un
ingr\'{e}dient essentiel de la preuve du th\'{e}or\`{e}me KMS est la
\emph{d\'{e}formation} des surfaces de translation,
par l'action du groupe $\Sl$ (et plus pr\'{e}cis\'{e}ment par son
sous-groupe constitu\'{e} des matrices diagonales). Ce processus de
d\'{e}formation est ais\'{e} \`{a} mettre en place dans le cadre des surfaces
de translation (voir la section~\ref{sec:action}), alors qu'il est
impossible dans l'espace des billards rationnels~! Moralement, cela
s'explique par le fait qu'un billard rationnel et son flot sont des
notions euclidiennes, alors qu'une surface de translation est une
notion affine. C'est ce changement de groupe structurel qui
permet de faire marcher la preuve. Ainsi, les surfaces de
translation sont un outil indispensable dans l'\'{e}tude des billards
rationnels.

\numero Nous terminons cette section par quelques remarques.
Dans l'\'{e}tude des surfaces de translation, il existe trois grandes
classes de r\'{e}sultats.
\begin{enumerate}
\item Tout d'abord, certains types de r\'{e}sultats sont valables pour
\emph{toutes} les surfaces de translation. C'est par exemple le cas
du th\'{e}or\`{e}me KMS. Les techniques utilis\'{e}es pour
d\'{e}montrer ce genre de r\'{e}sultats ne sont en fait pas tr\`{e}s
nombreuses (elles sont essentiellement g\'{e}om\'{e}triques et
combinatoires). La preuve du th\'{e}or\`{e}me KMS que nous allons
donner est assez repr\'{e}sentative de ce type d'arguments.
\item Nous verrons plus loin que l'on peut mettre une mesure sur l'ensemble des
surfaces de translation lorsque l'on fixe le type topologique (genre
et singularit\'{e}s). Certains r\'{e}sultats sont d\'{e}montr\'{e}s pour
\emph{presque toutes} les surfaces de translation, au sens de cette
mesure. Les r\'{e}sultats obtenus ainsi sont souvent beaucoup plus
pr\'{e}cis que ceux qui sont valables en toute
g\'{e}n\'{e}ralit\'{e}. N\'{e}anmoins, le prix \`{a} payer est que ces
r\'{e}sultats ne s'appliquent pas aux billards puisque l'ensemble
$\{\Su(P),\ P \textrm{ polygone rationnel} \}$ est de mesure nulle
dans l'espace des surfaces de translation~!
\item Enfin, on peut se restreindre \`{a} certaines classes de surfaces
(surfaces de Veech, ou surfaces de genre 2 par exemple) et chercher
\`{a} obtenir des r\'{e}sultats de classification, presque de nature
alg\'{e}brique.
\end{enumerate}

\noindent Les techniques utilis\'{e}es pour d\'{e}montrer des
r\'{e}sultats du deuxi\`{e}me ou du troisi\`{e}me type sont tr\`{e}s
vari\'{e}es, et ne seront pas du tout abord\'{e}es dans ce texte.

\section{Le th\'{e}or\`{e}me KMS faible}

\begin{Theorem}[KMS faible]
\label{theo:kms:faible}
Pour toute surface de translation $\Su$ et pour Lebesgue presque
tout $\theta\in\mathbb S^1$, le flot directionnel $\mathcal F_\theta$
sur $\Su$ est minimal.
\end{Theorem}

Ce th\'{e}or\`{e}me est d\^{u} \`{a} A.~Katok et A.~Zemljakov (\cite{KaZa}). Il
implique le corollaire suivant, qui est une version
``transitive'' du corollaire~\ref{cor:ergodique}.

\begin{Corollary}
L'ensemble $\{ Q\ ; \ Q \textrm{ est un polygone \`{a} } k \textrm{
  c\^{o}t\'{e}s et } \tilde \F \textrm{ est transitif} \}$ forme un
$G_\delta$ dense dans l'ensemble des polygones \`{a} $k$ c\^{o}t\'{e}s
(pour la topologie produit, c'est \`{a} dire en identifiant
l'ensemble des polygones \`{a} $k$ c\^{o}t\'{e}s avec un ouvert de
$\R^{2k}$).
\end{Corollary}

\begin{Remark}
Le r\'{e}sultat ci-dessus est bien plus faible que le th\'{e}or\`{e}me de
Kerckhoff, Masur et Smillie (th\'{e}or\`{e}me KMS fort). En
particulier, nous allons montrer ici que l'ensemble des
directions non-minimales est d\'{e}nombrable, donc de mesure nulle.
Dans le cas g\'{e}n\'{e}ral, il peut arriver que l'ensemble des
directions non uniquement ergodiques soit bien plus ``gros''.
Par exemple, pour les surfaces $\Su_\lambda$ de la
figure~\ref{fig:exemple:4:slit}, on peut montrer que si
$\lambda$ est diophantien alors l'ensemble des directions
non-uniquement ergodiques sur la surface $\Su_\lambda$ est de
dimension de Hausdorff $1/2$ (\cite{Cheung}, voir aussi
section~\ref{sec:raffinements}).
\end{Remark}

\noindent Le th\'{e}or\`{e}me~\ref{theo:kms:faible} repose sur les
propositions~\ref{prop:denombrable} et~\ref{prop:minimal} suivantes.

\begin{Proposition}
\label{prop:denombrable}
Sur une surface de translation, l'ensemble des directions de
connexions de selles est au plus d\'{e}nombrable.
\end{Proposition}

\begin{proof}
Il suffit de passer au rev\^{e}tement universel.
\end{proof}

\begin{Proposition}
\label{prop:minimal}
Soit $\Su$ une surface de translation. Si le flot directionnel
$\mathcal F_\theta$ ne poss\`{e}de pas de connexion de selles alors il
est minimal~: c'est \`{a} dire toutes les feuilles de $\mathcal
F_\theta$ sont denses.
\end{Proposition}

\noindent On peut sans perte de g\'{e}n\'{e}ralit\'{e} se restreindre au
cas du feuilletage vertical. Avant de d\'{e}montrer cette
proposition, nous aurons besoin d'un lemme g\'{e}om\'{e}trique et d'un
lemme technique.

\begin{Lemma}
\label{lm:geo}
Soit $\Su$ une surface de translation de genre $g\geq 2$. Supposons
que le flot vertical poss\`{e}de une feuille verticale r\'{e}guli\`{e}re
ferm\'{e}e. Alors le flot vertical poss\`{e}de aussi une connexion de
selles.
\end{Lemma}

\begin{proof}
Soit $\alpha$ une feuille verticale r\'{e}guli\`{e}re ferm\'{e}e sur $\Su$.
Comme $\Su$ ne poss\`{e}de qu'un nombre fini de singularit\'{e}s, il
existe une autre feuille verticale ferm\'{e}e parall\`{e}le \`{a} $\alpha $
et proche de $\alpha$. De proche en proche, cela fournit un
plongement isom\'{e}trique du cylindre $]0,r[ \times \mathbb S^1$
de p\'{e}rim\`{e}tre $|\alpha|$ et de largeur $r$ dans $\Su$. Ce
cylindre est feuillet\'{e} par des feuilles r\'{e}guli\`{e}res verticales
homologues \`{a} $\alpha$. Par ailleurs $\Aire(\Su)=1$ donc on ne
peut pas ``\'{e}paissir'' ce cylindre \`{a} l'infini. Les deux seules
obstructions topologiques sont les suivantes. Ou bien il existe
une connexion de selles verticale, situ\'{e}e sur le bord du
cylindre, ou bien $\Su$ s'obtient en recollant les deux bords
du cylindre ensemble. Ce dernier cas implique alors que $g(\Su)
= 1$ ce qui est une contradiction. Le lemme est ainsi d\'{e}montr\'{e}.
\end{proof}

\begin{Lemma}
\label{lm:technique}
Soit $\alpha^{+}$ une feuille verticale (orient\'{e}e positivement)
non p\'{e}riodique (i.e. non ferm\'{e}e). Soient $P\in \alpha$ un point
et $I=[P,Q]$ un segment g\'{e}od\'{e}sique transverse au flot vertical.
Alors $]P,Q[ \ \cap\ \alpha^{+} \not = \emptyset$.
\end{Lemma}

\begin{proof}
La surface $\Su$ ne poss\`{e}de qu'un nombre fini de s\'{e}paratrices
verticales. Ainsi on peut choisir $I'=[P,Q'] \subset I$ tel que
toutes les feuilles verticales issues de $I'$ ne rencontrent
pas de singularit\'{e}s avant de revenir \`{a} $I$. Consid\'{e}rons alors
un petit rectangle vertical $R_h$ {\it plong\'{e}} dans $\Su$, de
base $I'$ et de hauteur $h$, avec $h$ petit. Le bord vertical
gauche de $R_h$ est inclus dans $\alpha^{+}$. Notons par
$\beta^{+}$ le bord vertical droit. Lorsque $h \to \infty$ le
rectangle $R_h$ ne reste pas plong\'{e} dans $\Su$~: en effet
l'aire de $\Su$ est finie et l'aire de $R_h$ tends vers
l'infini. Donc (rappelons que $R_h$ ne rencontre pas de
singularit\'{e} avant de revenir \`{a} $I$) deux possibilit\'{e}s peuvent
arriver. Ou bien $\alpha^{+}$ intersecte $]P,Q'[$ et le lemme
est d\'{e}montr\'{e} (figure~\ref{fig:lemme}a), ou bien c'est
$\beta^{+}$ qui intersecte $]P,Q'[$ (figure~\ref{fig:lemme}b).

\begin{figure}[htbp]
\begin{center}

\psfrag{r}{$\scriptstyle R_h$} \psfrag{a}{$\scriptstyle
\alpha^{+}$} \psfrag{b}{$\scriptstyle \beta^{+}$}
\psfrag{p}{$\scriptstyle P$} \psfrag{q}{$\scriptstyle Q'$}
\psfrag{qp}{$\scriptstyle Q$} \psfrag{qpp}{$\scriptstyle Q''$}

  \subfigure[]{\epsfig{figure=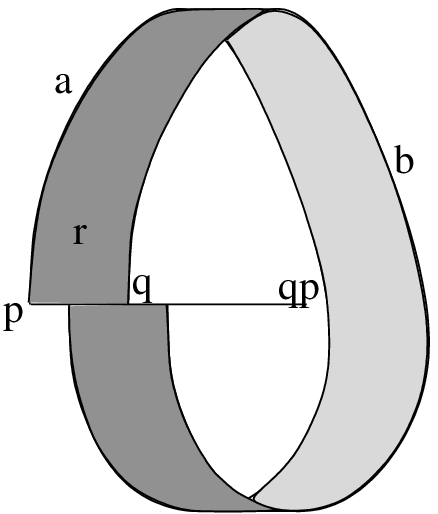,width=5.0cm}} \qquad \qquad
  \subfigure[]{\epsfig{figure=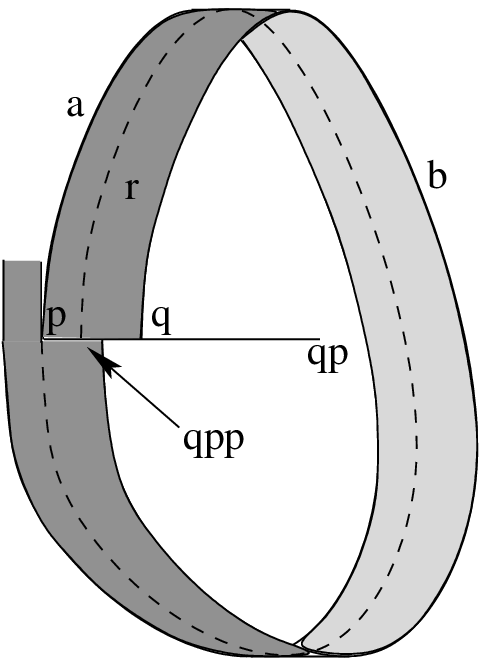,width=5.0cm}}

\end{center}
\caption{\label{fig:lemme}
}
\end{figure}

\noindent Dans ce dernier cas, notons $Q'' \in I'$ l'unique
point tel que la feuille verticale issue de $Q''$ dans $R_h$
arrive \`{a} $P$. Consid\'{e}rons alors $I''=[P,Q''] \subset I'$ et
appliquons le raisonnement analogue \`{a} celui pr\'{e}c\'{e}dent. On
obtient alors un petit rectangle de base $I''$. Le m\^{e}me
argument am\`{e}ne \`{a} la dichotomie pr\'{e}c\'{e}dente. Il est facile de
voir alors que la deuxi\`{e}me conclusion dans cette dichotomie ne
peut-\^{e}tre satisfaite. Le lemme est d\'{e}montr\'{e}.
\end{proof}

\numero On peut maintenant conclure.
\begin{proof}[Preuve de la proposition~\ref{prop:minimal}]
Supposons que $\Su$ ne contienne pas de connexion de selles
verticale. Soit $\beta$ une feuille verticale. Si $\beta$ est
ferm\'{e}e alors on obtient une contradiction gr\^{a}ce au
lemme~\ref{lm:geo}. Notons $A$ l'adh\'{e}rence de $\beta$ dans
$\Su$. Nous allons montrer que $A=\Su$. Sinon prenons $P \in
\partial A$. Soit $\alpha$ la feuille verticale issue de $P$.
Si $\alpha$ est ferm\'{e}e alors le lemme~\ref{lm:geo} fournit de
nouveau une contradiction. Sinon on note $I=[P,Q]$ un petit
segment g\'{e}od\'{e}sique, transverse au flot vertical, passant par
$P$ et tel que son int\'{e}rieur $\Int(I)$ soit inclus dans $\Su
\setminus A$. Notons aussi que $A$ est invariant par le flot
vertical donc $\alpha \subset A$. Le lemme~\ref{lm:technique}
implique alors que $\Int(I) \cap \alpha \not =\emptyset$ ce qui
est \'{e}videmment incompatible avec les deux derni\`{e}res assertions.
Donc $A=\Su$ et ainsi toutes les feuilles verticales sont
denses. La proposition est d\'{e}montr\'{e}e.
\end{proof}

\begin{Remark}
Pour aller plus loin que le th\'{e}or\`{e}me~\ref{theo:kms:faible}, on peut
chercher \`{a} caract\'{e}riser les directions exceptionnelles,
c'est-\`{a}-dire les directions pour lesquelles il y a une connexion
de selles. C'est un probl\`{e}me tr\`{e}s difficile en g\'{e}n\'{e}ral. Par
exemple en genre $2$, pour les surfaces de Veech, cet ensemble est
toujours (quitte \`{a} normaliser convenablement)
un corps de nombres union l'infini (plus pr\'ecis\'ement, 
soit $\mathbb Q \cup \{\infty \}$, soit $K \cup \{\infty \}$ o\'u $K$ 
est un corps quadratique r\'{e}el).
\end{Remark}

\section{Preuve du th\'{e}or\`{e}me KMS}

\begin{Definition}
\label{def:PartirInfini} Soit $\Su$ une surface de translation.
On dira que $\Su$ part \`{a} l'infini dans la direction $\pi/2$, ou
simplement que $\Su$ part \`{a} l'infini, si pour tout
$\epsilon>0$, il existe $T>0$ tel que, pour tout $t\geq T$, il
existe une connexion de selles dans $\Su$ (de composantes
horizontales et verticales $h$ et $v$) telle que $e^t |h| +
e^{-t}|v| \leq \epsilon$.

\noindent On parlera de m\^{e}me de ``partir \`{a} l'infini dans la direction
$\theta$'' en consid\'{e}rant les composantes des connexions de selle
dans la direction $\theta$ et la direction orthogonale \`{a} $\theta$.
\end{Definition}

\numero Les deux ingr\'{e}dients essentiels de la preuve du th\'{e}or\`{e}me
KMS sont les r\'{e}sultats suivants.

\begin{ThmA}
Soit $\Su$ une surface de translation qui ne part pas \`{a} l'infini
dans la direction $\theta$. Alors le flot $\F_\theta$ est uniquement
ergodique.
\end{ThmA}

\begin{ThmB}
Soit $\Su$ une surface de translation. Alors l'ensemble
  \begin{equation*}
  \{\theta\in \Sbb^1,\ \Su \textrm{ part \`{a} l'infini
  dans la direction }\theta\}
  \end{equation*}
est de mesure nulle.
\end{ThmB}

\noindent Ces deux th\'{e}or\`{e}mes impliquent manifestement le th\'{e}or\`{e}me
KMS. Leurs d\'{e}monstrations sont compl\`{e}tement diff\'{e}rentes et
ind\'{e}pendantes, et seront pr\'{e}sent\'{e}es dans les
sections~\ref{PreuveThmA} et~\ref{PreuveThmB} respectivement.
Remarquons aussi tout de suite que le th\'{e}or\`{e}me~$A$
n'est pas une \'{e}quivalence, contrairement au cas du genre~$1$ (voir
section~\ref{sec:raffinements}). Le th\'{e}or\`{e}me~$A$ est aussi
connu sous le nom de crit\`{e}re de Masur. \medskip

Notons qu'il existe une version combinatoire de ce crit\`{e}re, due
\`{a} M.~Boshernitzan (voir~\cite{Boshernitzan}). Nous reviendrons
section~\ref{sec:raffinements} sur ce crit\`{e}re. \medskip

Avant de d\'{e}montrer les th\'{e}or\`{e}mes $A$ et $B$, on va
introduire l'action de $\Sl$ sur les surfaces de translation, ce qui
permet de reformuler de mani\`{e}re un peu plus confortable la notion de
d\'{e}part \`{a} l'infini d\'{e}finie en~\ref{def:PartirInfini}.

\section{Action de $\Sl$ sur les surfaces de translation}
\label{sec:action}

\numero Le groupe $\Sl$ agit lin\'{e}airement sur
$\R^2$~; il agit donc aussi naturellement sur les surfaces de
translation via les cartes locales. Si $(\Su,\Sigma,\omega)$ est une
surface de translation avec $\omega=\left\{ (U_i,z_i) \right\}_i$ et
si $A\in \Sl$ est une matrice alors on d\'{e}finit l'action comme suit.
  \begin{equation}
  A \cdot (\Su,\Sigma,\omega) := (\Su,\Sigma,A\omega)
  \end{equation}
o\`{u} $A\omega$ est, par d\'{e}finition, le nouvel atlas plat $A
\omega:=\left\{ (U_i, A \circ z_i) \right\}_i$, $A \circ z_i$
d\'{e}signant l'action lin\'{e}aire de $A$ sur la coordonn\'{e}e $z_i$.
Remarquons que l'action pr\'{e}serve les surfaces d'aire $1$.

\numero Nous serons particuli\`{e}rement int\'{e}ress\'{e} par les $3$
sous-groupes \`{a} $1$ param\`{e}tre suivants.
  $$
  g_t = \left( \begin{array}{ll}
  e^t & 0 \\
  0  & e^{-t}
  \end{array} \right), \qquad
  h_s = \left( \begin{array}{ll}
  1 & s \\
  0  & 1
  \end{array} \right) \qquad \textrm{ et } \qquad
  R_\theta = \left( \begin{array}{lr}
  \cos(\theta) & -\sin(\theta) \\
  \sin(\theta)  & \cos(\theta)
  \end{array} \right).
  $$
L'action de $g_t$ sera
appel\'{e}e {\it flot g\'{e}od\'{e}sique de Teichm\"{u}ller}
ou simplement flot g\'{e}od\'{e}sique pour une raison qui deviendra claire
dans les sections suivantes. L'action de $g_t$ ``\'{e}crase'' les
feuilles verticales par $e^t$ et ``dilate'' les feuilles
horizontales par $e^t$. Le flot $h_s$ est le flot horocyclique et le
flot $R_\theta$ est le flot circulaire. L'action de $R_\theta$
consiste \`{a} faire tourner le feuilletage vertical et le feuilletage
horizontal de $\Su$ d'un angle de $-\theta$.

\numero Nous noterons pour la suite $\Su_\theta = R_\theta \Su$
et $\Su_{\theta,t} = g_t \Su_\theta$. On v\'{e}rifie sans peine que
le flot vertical sur $\Su_{\pi/2-\theta}$ correspond au flot
dans la direction $\theta$ sur $\Su$.

\begin{Definition}
Si $\Su$ est une surface de translation, on appelle \emph{systole} de
$\Su$ la longueur (dans la m\'etrique euclidienne) 
de la plus petite connexion de selles de $\Su$. On la notera $\sys(\Su)$.
\end{Definition}

\numero Avec ces d\'{e}finitions, on peut reformuler ais\'{e}ment la
d\'{e}finition~\ref{def:PartirInfini} comme suit~: une surface part
\`{a} l'infini si et seulement si $\sys(g_t \Su) \to 0$ quand $t
\to \infty$. Elle part \`{a} l'infini dans la direction $\theta$ si
et seulement si $\sys(g_t R_{\pi/2-\theta} \Su)\to 0$ quand $t
\to \infty$.

\section{Preuve du th\'{e}or\`{e}me A}
\label{PreuveThmA}

Comme d'habitude, il suffit de le prouver pour le flot vertical
$\theta=\pi/2$, puisque le r\'{e}sultat g\'{e}n\'{e}ral s'en d\'{e}duit par
rotation. L'id\'{e}e essentielle de la preuve est la suivante~: si la
surface ne part pas \`{a} l'infini, alors on va pouvoir extraire une
``sous-suite convergente'' (en un certain sens \`{a} expliquer) par un
argument de compacit\'{e}. On conclura ensuite en utilisant la
connexit\'{e} de la surface limite.

\begin{Definition}
\label{def:Convergence}
Soient $s,g\in \N$. On dira qu'une suite de surfaces de translation
$\Su_n$ d'aire~$1$ (de genre $g$, avec $s$ singularit\'{e}s) converge
s'il existe des triangulations de $\Su_n$ par des triangles
$T_n^1,\dots, T_n^{f}$ (avec $f$ ne d\'{e}pendant que de $g$ et $s$,
d'apr\`{e}s~\eqref{FormuleFaces}) telles que
\begin{enumerate}
\item chaque suite de triangle $(T_n^i)_{n\in \N}$ converge, comme suite de
triangles dans $\R^2$, vers un triangle (peut-\^{e}tre d\'{e}g\'{e}n\'{e}r\'{e}
i.e. aplati) $T^i$.
\item le motif de recollement des $T_n^i$ pour former $\Su_n$ est
constant \`{a} partir d'un certain rang.
\end{enumerate}
On impose aussi dans ce type de convergence que la systole $\sys(\Su_n)$ de
$\Su_n$ soit uniform\'{e}ment minor\'ee par une constante non-nulle. Ceci implique
en particulier que la surface limite est non ``d\'{e}g\'{e}n\'{e}r\'{e}e''.
\end{Definition}

\noindent Il est alors possible de recoller les triangles $T^i$
suivant les m\^{e}mes motifs que les $T_n^i$ pour $n$ assez grand, pour
obtenir une surface de translation, que l'on notera $\Su_\infty$, de
genre $g$ avec $s$ singularit\'{e}s. On dira alors que $\Su_n$
converge vers $\Su_\infty$. La topologie de $\Su_\infty$ est la m\^{e}me
que celle des surfaces $\Su_n$ (elle est d\'{e}termin\'{e}e par le motif
combinatoire de recollement des triangles). \medskip

\noindent Notons aussi que si l'on n'imposait pas la condition sur la
systole de $\Su_n$, il serait \emph{a priori} possible que des
singularit\'{e}s fusionnent dans la surface limite~! Par exemple la 
surface limite pourrait ne plus \^etre connexe.

\numero
\label{par:PossibilitesConvergence}
Ce type de convergence est assez faible, mais il permet de r\'{e}aliser
les op\'{e}rations suivantes.
\begin{enumerate}
\item Si $x_n \in \Su_n$, cela a un sens de parler de la convergence
de $x_n$ vers $x\in \Su_\infty$, en regardant ce qui se passe dans le
triangle contenant $x_n$.
\item Si on se donne une partie de $\Su_\infty$, on peut l'approcher
par une partie de $\Su_n$ si $n$ est assez grand.
\end{enumerate}
\noindent Nous n'aurons pas besoin d'autre chose pour d\'{e}montrer le
th\'{e}or\`{e}me~$A$.

\begin{Lemma}
\label{lem:BonneTriangulation}
Soit $\epsilon_0>0$. Il existe une constante $C_0$ ne d\'{e}pendant que
de $\epsilon_0$ telle que toute surface de translation
d'aire $1$ et de systole au moins $\epsilon_0$ admette une
triangulation dont tous les c\^{o}t\'{e}s ont une longueur au plus $C_0$.
\end{Lemma}

\noindent Par triangulation, on entend ici une triangulation dont les
sommets sont les singularit\'{e}s et les ar\^{e}tes des connexions de
selles, comme en~\ref{def:Triangulation}.

\begin{proof}
Toute surface de translation admet une triangulation,
d'apr\`{e}s~\ref{ExisteTriangulation}. Partons d'une telle
triangulation, et consid\'{e}rons une ar\^{e}te de longueur $K$ maximale.
Notons-la $AC$, et consid\'{e}rons les deux triangles qui
bordent cette ar\^{e}te, comme sur la figure \ref{triangle}.
\begin{figure}[htbp]

\begin{center}
\psfrag{a}{$A$} \psfrag{b}{$B$}
\psfrag{c}{$C$} \psfrag{bp}{$B'$}
\psfrag{h}{$H$} \psfrag{hp}{$H'$}

\includegraphics[width=8.5cm]{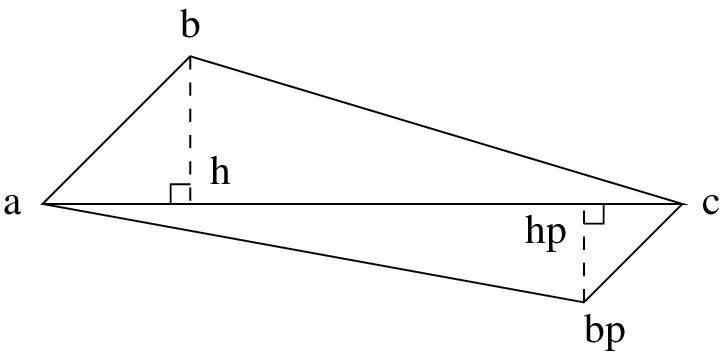}
\end{center}
\caption{\label{triangle}}
\end{figure}

\noindent Comme l'aire de la surface est au plus $1$, ces deux
triangles sont d'aire au plus $1$. Ainsi l'aire de $AHB$ est
$1/2(AC.BH) \leq 1$ ce qui donne $BH\leq 2/K$. Comme $AB$ et
$BC$ sont de longueur au moins $\epsilon_0$, Pythagore dans le
triangle $ABH$ donne $AH^2 = AB^2-BH^2 \geq \epsilon^2_0 -
4/K^2 \geq \epsilon^2_0/4$ si $K$ est grand (c'est \`{a} dire si $K
\geq \cfrac{4}{\sqrt{3}\epsilon_0}$). Donc $AH\geq
\epsilon_0/2$. De la m\^{e}me mani\`{e}re on tire $CH\geq
\epsilon_0/2$. Les m\^{e}mes in\'{e}galit\'{e}s ont lieu pour $B'H'$, $AH'$
et $CH'$. \medskip

On retriangule alors la surface en retirant
l'ar\^{e}te $AC$ et en la rempla\c{c}ant par une ar\^{e}te $BB'$. Cette
ar\^{e}te est de longueur au plus $BH+B'H'+HH' \leq 4/K +(K-\epsilon_0)$,
qui est major\'{e} par $K-\epsilon_0/2$ si $K$ est tr\`{e}s grand (c'est \`{a} dire si
$K \geq 2/\epsilon_0$). \medskip

Ainsi, on peut faire d\'{e}cro\^{\i}tre la longueur maximale des
ar\^{e}tes de la triangulation d'une quantit\'{e} fix\'{e}e. En r\'{e}p\'{e}tant ce
processus, on arrive alors \`{a} une triangulation satisfaisant la borne requise.
\end{proof}

\begin{Corollary}
\label{cor:Compacite}
Soit $\Su_n$ une suite de surfaces de translation d'aire $1$, avec un
genre et un nombre de singularit\'{e}s fix\'{e}s. On suppose que la
systole de $\Su_n$ ne tend pas vers $0$ quand $n$ tend vers
l'infini. Alors $\Su_n$ admet une sous-suite qui converge vers
$\Su_\infty$ au sens de la d\'{e}finition~\ref{def:Convergence}. De plus
le genre et le nombre de singularit\'{e}s de $\Su_\infty$ sont les
m\^{e}mes que $\Su_n$.
\end{Corollary}

\begin{proof}
Gr\^{a}ce au lemme~\ref{lem:BonneTriangulation}, on obtient une
sous-suite de $\Su_n$ ayant une triangulation avec des c\^{o}t\'{e}s de
taille born\'{e}e sup\'{e}rieurement et inf\'{e}rieurement, par des
triangles $T_n^i$. Les suites de tels triangles dans $\R^2$ forment un
ensemble compact, on peut donc choisir une sous-suite de $\Su_n$ de
telle sorte que les $T_n^i$ convergent. Il y a aussi un nombre fini de
motifs de recollement possibles, on obtient donc une sous-suite
convergente quitte \`{a} extraire encore une fois. Notons aussi que dans
le lemme
pr\'{e}c\'{e}dent, lorsque l'on retriangule la surface, la nombre de faces
$f$ et le genre $g$ de la surface restent inchang\'{e}s. Donc d'apr\`{e}s
la formule~\eqref{FormuleFaces} le nombre de singularit\'{e}s $s$ reste
aussi inchang\'{e}. Le corollaire est d\'{e}montr\'{e}.
\end{proof}

\numero On peut commencer la preuve du th\'{e}or\`{e}me~$A$ proprement
dite. \medskip

\noindent Soit $\Su$ une surface de translation. On note
$\Su_t$ la surface obtenue en dilatant les distances
horizontales d'un facteur $e^t$, et en contractant les
distances verticales de ce m\^{e}me facteur. Formellement, $\Su_t=
g_t \Su$. On suppose que $\Su$ ne part pas \`{a} l'infini. Ceci
garantit alors l'existence d'une suite $t_n\to \infty$ et de
$\epsilon_0>0$ tels que toute connexion de selles de la surface
$\Su_{t_n}$ soit de longueur au moins $\epsilon_0$. D'apr\`{e}s le
corollaire~\ref{cor:Compacite}, on peut extraire de la suite
$\Su_{t_n}$ une sous-suite convergente, que nous noterons
simplement $\Su_n$.
\medskip

\noindent Les surfaces $\Su_n$ et $\Su$ sont des structures de
translation sur la m\^{e}me surface topologique sous-jacente. Pour
bien les distinguer, nous noterons $f_n$ l'application $f_n :
\Su \rightarrow \Su_n$ qui est l'identit\'{e} sur cette surface
topologique sous-jacente. \medskip

\noindent Soit $\mu$ une mesure invariante ergodique pour le
flot vertical sur $\Su$. On note $B(\mu) \subset \Su$
l'ensemble des points $x$ qui sont typiques pour la mesure
$\mu$, i.e., tels que $\frac{1}{T} \int_{0}^T
\delta_{\F_{\pi/2}^t(x)}$ converge faiblement vers $\mu$ quand
$T\to \pm \infty$, o\`{u} $\delta_y$ d\'{e}signe la masse de Dirac au
point $y$.

\noindent On d\'{e}finit un sous-ensemble $A(\mu) \subset \Su_\infty$
comme suit. C'est la r\'{e}union, pour $x\in B(\mu)$, de l'ensemble des
valeurs d'adh\'{e}rence possibles dans $\Su_\infty$ de la suite
$f_n(x)\in \Su_n$ (comme d\'{e}fini
en~\ref{par:PossibilitesConvergence}). \medskip

\noindent Pour d\'{e}montrer qu'il y a une seule mesure invariante
$\mu$, on va en fait montrer qu'il y a un seul ensemble $A(\mu)$. On
appellera \emph{rectangle} un rectangle dont les c\^{o}t\'{e}s sont horizontaux et
verticaux, et \emph{rectangle plong\'{e}} un rectangle qui ne rencontre
pas les singularit\'{e}s.

\begin{Lemma}
\label{lem:MuEgaleNu}
Consid\'{e}rons dans $\Su_\infty$ un rectangle plong\'{e}. On suppose
qu'il existe $x\in B(\mu)$ et $y\in B(\nu)$ (o\`{u} $\mu$ et $\nu$ sont
deux mesures invariantes) tels que $f_n(x)$ et $f_n(y)$ convergent
respectivement vers deux sommets oppos\'{e}s de ce rectangle. Alors
$\mu=\nu$.
\end{Lemma}

\begin{proof}
Supposons que $\mu$ et $\nu$ ne co\"{\i}ncident pas. Comme les
rectangles plong\'{e}s de hauteur $1$ engendrent la tribu des
bor\'{e}liens, il existe un rectangle $U=I\times J$, plong\'{e} dans
$\Su$, tel que $\mu(U) \not=\nu(U)$, et de hauteur $1$.
\medskip

Si $n$ est assez grand, on peut tracer dans $\Su_n$ un
rectangle $R_n$ \`{a} c\^{o}t\'{e}s horizontaux et verticaux, ne contenant
pas de singularit\'{e}, et dont $f_n(x)$ et $f_n(y)$ sont deux sommets
oppos\'{e}s. Cela r\'{e}sulte de l'hypoth\`{e}se dans $\Su_\infty$, et de la
convergence de $\Su_n$ vers $\Su_\infty$. De plus, la taille verticale
$c_n$ de $R_n$ converge vers $c>0$. \medskip

\noindent Le rectangle $f_n^{-1}(R_n)$ dans $\Su$ est un
rectangle de hauteur $T_n=c_n e^{t_n}\to \infty$, il est tr\`{e}s
haut et tr\`{e}s fin, et poss\`ede $x$ et $y$ pour sommets oppos\'{e}s. Notons
$C_n(x)$ et $C_n(y)$ les c\^{o}t\'{e}s verticaux de $f_n^{-1}(R_n)$ contenant
respectivement $x$ et $y$. Comme $x$ est typique pour $\mu$, le
th\'{e}or\`{e}me de Birkhoff donne $\Card( I \cap C_n(x)) \sim T_n
\mu(U)$ quand $n\rightarrow \infty$. 
De m\^{e}me, $\Card( I \cap C_n(y)) \sim T_n \nu(U)$ quand $n\rightarrow \infty$.
Comme le rectangle $f_n^{-1}(R_n)$ ne contient pas de
singularit\'{e}, tout segment horizontal traversant $C_n(x)$ doit
traverser aussi $C_n(y)$ (sauf \'{e}ventuellement une fois pr\`{e}s de
chaque bord du segment horizontal, voir
figure~\ref{fig:rectangle}).

\begin{figure}[htbp]
\begin{center}
\psfrag{Cx}{$\scriptstyle C_n(x)$}  \psfrag{Cy}{$\scriptstyle C_n(y)$}
\psfrag{x}{$\scriptstyle x$}  \psfrag{y}{$\scriptstyle y$}
\psfrag{i}{$\scriptstyle I$} \psfrag{j}{$\scriptstyle J$}
\psfrag{fn}{$\scriptstyle f_n^{-1}(R_n)$}

\includegraphics[width=6.5cm]{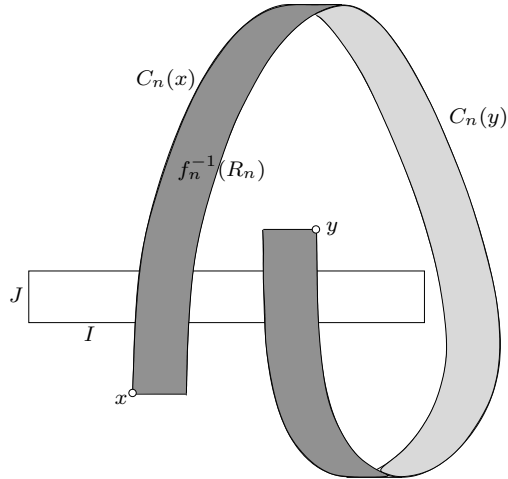}
\end{center}
\caption{ \label{fig:rectangle} Le rectangle $f_n^{-1}(R_n)$ vu
sur la surface $\Su$. Un segment horizontal traversant $C_n(x)$
traverse aussi $C_n(y)$ (sauf \'{e}ventuellement une fois pr\`{e}s de
chaque bord du segment horizontal) donc $|\Card(I\cap C_n(x))-
\Card( I\cap C_n(y))| \leq 2$. }
\end{figure}

\noindent Ceci implique alors que $|\Card(I\cap C_n(x))- 
\Card( I\cap C_n(y))| \leq 2$. Comme $T_n \to \infty$, 
on obtient $T_n \mu(U) \sim T_n \nu(U)$ et donc 
$\mu(U)=\nu(U)$, ce qui est absurde et conclut.
\end{proof}

\begin{Lemma}
\label{lem:AmuPartout}
Tout ouvert de $\Su_\infty$ contient un point appartenant \`{a} un
ensemble $A(\mu)$ pour une certaine mesure invariante ergodique $\mu$.
\end{Lemma}

\begin{proof}
Soit $U$ un ouvert de $\Su_\infty$~; c'est la limite d'une suite
d'ouverts $U_n \in \Su_n$. De plus, $\Leb(U_n) \to \Leb(U)>0$. Notons
$V_n=f_n^{-1}(U_n) \subset \Su$. Comme l'application $f_n : \Su \to
\Su_n$ envoie la mesure de Lebesgue de $\Su$ sur la mesure de
Lebesgue de $\Su_n$, $\Leb(V_n)$ ne tend pas vers $0$. \medskip

\noindent Montrons qu'il existe une mesure $\mu$ invariante
ergodique telle que $\mu(V_n)$
ne tende pas vers $0$, par l'absurde. On \'{e}crit la d\'{e}sint\'{e}gration
ergodique de $\Leb$, comme $\Leb= \int_{\Omega} \mu_\xi \D P(\xi)$
o\`{u} $(\Omega, P)$ est un espace probabilis\'{e} et $(\mu_\xi)_{\xi\in
\Omega}$ parcourt l'ensemble des mesures invariantes ergodiques. Si
tous les $\mu_\xi(V_n)$ tendaient vers $0$, alors par convergence
domin\'{e}e $\Leb(V_n)$ tendrait \'{e}galement vers $0$, ce qui est
absurde. \medskip

\noindent Fixons donc $\mu$ telle que $\mu(V_n)$ ne tende pas vers
$0$. Alors
 \begin{equation*}
  \mu\{ x\tq x\text{ appartient \`{a} une infinit\'{e} de }V_n\}>0.
  \end{equation*}
En particulier, il existe $x\in B(\mu)$ qui appartient \`{a} une
infinit\'{e} de $V_n$. La suite $f_n(x) \in U_n$ admet alors une valeur
d'adh\'{e}rence dans $\overline{U}$.
\end{proof}

\numero On peut maintenant conclure.

\begin{proof}[Preuve du th\'{e}or\`{e}me~$A$]
On fixe des petits ouverts $U_1,\dots, U_N$ de la surface
limite $\Su_\infty$, satisfaisant la propri\'{e}t\'{e} suivante~: pour
tous points $x_i \in U_i$~:
\begin{enumerate}
\item pour tous $i,j$, il existe une suite de rectangles plong\'{e}s
ayant des sommets oppos\'{e}s parmi les $x_k$,
tels que le premier rectangle ait $x_i$ pour sommet et le dernier ait
$x_j$ pour sommet.
\item pour tout $x\in \Su_\infty$ qui n'est pas une singularit\'{e}, il
existe un rectangle plong\'{e} ayant $x$ et l'un des $x_j$ pour sommets
oppos\'{e}s.
\end{enumerate}

\noindent De tels ouverts existent bien, par un argument de
connexit\'{e} sur la surface $\Su_\infty$. L'ouvert $U_1$ contient un
point $x_1$ appartenant \`{a} un ensemble $A(\mu_1)$,
par le lemme~\ref{lem:AmuPartout}. Quitte \`{a} consid\'{e}rer une
sous-suite de $\Su_n$, on peut m\^{e}me supposer qu'il existe $y_1 \in
B(\mu)$ tel que $f_n(y_1)$ converge vers $x_1$. On se restreint \`{a}
cette sous-suite d'indices. L'ouvert $U_2$ contient un point $x_2$
appartenant \`{a} un ensemble $A(\mu_2)$, encore par le
lemme~\ref{lem:AmuPartout}, et quitte \`{a} extraire on peut encore
supposer
qu'on a convergence le long de la suite $\Su_n$. 
On r\'{e}p\`{e}te ainsi ce processus d'extraction
$N$ fois, et on obtient $N$ points $x_1,\dots,x_N$ dans $U_1,\dots,
U_N$, limites de suites $f_n(y_i)$ avec $y_i \in B(\mu_i)$ pour une
certaine mesure invariante ergodique $\mu_i$. \medskip

\noindent La premi\`{e}re propri\'{e}t\'{e} dans le choix des ouverts $U_i$,
combin\'{e}e avec le lemme~\ref{lem:MuEgaleNu}, montre que toutes les
mesures $\mu_i$ sont \'{e}gales \`{a} une m\^{e}me mesure $\mu$. \medskip

\noindent Soit maintenant $\nu$ une autre mesure ergodique
invariante. Soit $y\in B(\mu)$, consid\'{e}rons $x$ une valeur d'adh\'{e}rence dans
$\Su_\infty$ de la suite $f_n(y)$, disons que $x$ est la limite d'une
suite $f_{j(n)}(y)$. Si $x$ n'est pas une singularit\'{e},
il existe un rectangle plong\'{e} dans $\Su_\infty$ ayant $x$ et l'un
des $x_i$ comme sommets. En appliquant le lemme~\ref{lem:MuEgaleNu}
\`{a} la suite $\Su_{j(n)}$, on obtient $\mu=\nu$. Si $x$ est une
singularit\'{e}, on ne peut pas appliquer le lemme~\ref{lem:MuEgaleNu},
mais on peut encore conclure en reprenant sa preuve et en utilisant le
fait que $f_{j(n)}(y)$ n'est pas une singularit\'{e} (les d\'{e}tails
g\'{e}om\'{e}triques sont laiss\'{e}s au lecteur). Dans tous les cas, on
obtient $\mu=\nu$. On a ainsi montr\'{e} qu'il y avait une seule mesure
invariante, ce qui conclut donc la preuve du th\'{e}or\`{e}me~$A$.
\end{proof}

\section{Preuve du th\'{e}or\`{e}me B}
\label{PreuveThmB}

Rappelons que nous devons d\'{e}montrer que, si $\Su$ est une surface de
translation, alors l'ensemble $\{\theta\in \Sbb^1,\ \Su \textrm{ part \`{a} l'infini
dans la direction }\theta\}$ est de mesure nulle.

\numero On va raisonner par l'absurde pour d\'{e}montrer le
th\'{e}or\`{e}me B. Donnons tout d'abord une traduction ``na\"{\i}ve'' de la
n\'{e}gation du th\'{e}or\`{e}me B.

Si ce th\'{e}or\`{e}me est faux pour une surface $\Su$, alors il existe
un sous-ensemble $S$ de $\Sbb^1$ de mesure de Lebesgue non
nulle tel que la surface $\Su$ tend vers l'infini dans la
direction $\theta$, pour tout $\theta \in S$. Autrement dit,
pour tout $\epsilon>0$ et pour tout $\theta \in S$, il existe
$T(\epsilon,\theta)$ tel que, pour tout $t>T(\epsilon,\theta)$,
la surface $\Su_{\theta,t}$ admet une connexion de selles de
longueur au plus $\epsilon$. On choisit $T(\epsilon,\theta)$
minimal, de telle sorte que la fonction $\theta \mapsto
T(\epsilon,\theta)$ est mesurable. Par cons\'{e}quent, il existe un
ensemble $S_\epsilon \subset S$ de mesure au moins $\Leb(S)/2$
tel que $T_\epsilon:=\sup_{\theta\in S_\epsilon}
T(\epsilon,\theta)$ soit fini. Toutes les surfaces
$\Su_{\theta, T_\epsilon+1}$ pour $\theta \in S_\epsilon$ ont
alors une connexion de selles de longueur au plus $\epsilon$.
En appliquant un proc\'{e}d\'{e} diagonal \`{a} la suite $\epsilon_n=1/n$,
on obtient la conclusion suivante.

\begin{Proposition}
\label{prop:RencontrePTpetit}
Soit $\Su$ une surface de translation qui ne satisfait pas la
conclusion du th\'{e}or\`{e}me B. Alors il existe des sous-ensembles $S_n$
de $\Sbb^1$ avec $\inf_{n\in\N} \Leb(S_n)>0$,
une suite $T_n \to \infty$ et une suite $\epsilon_n \to 0$ tels que :
pour tout $n\in \N$, pour tout $\theta \in S_n$, la surface
$\Su_{\theta, T_n}$ a une connexion de selles de longueur au plus
$\epsilon_n$.
\end{Proposition}

Malheureusement, cette n\'{e}gation du th\'{e}or\`{e}me B est trop na\"{\i}ve
pour amener directement \`{a} une contradiction. On aura besoin de la
version forte suivante~:

\begin{Proposition}
\label{prop:RencontrePasPetit}
Soit $\Su$ une surface de translation ne satisfaisant pas la
conclusion du th\'{e}or\`{e}me B. Alors il existe $c>0$, des sous-ensembles $S_n$
de $\Sbb^1$ avec $\inf_{n\in \N} \Leb(S_n)>0$,
une suite $T_n \to \infty$ et une suite $\epsilon_n \to 0$ tels que :
pour tout $n\in \N$, pour tout $\theta \in S_n$, la surface
$\Su_{\theta, T_n}$ a une connexion de selles de longueur au plus
$\epsilon_n$, qui ne rencontre pas d'autre connexion de selles de longueur
$\leq c$.
\end{Proposition}

Pour d\'{e}montrer cette proposition, on devra travailler avec des
ensembles de connexions de selles courtes et deux \`{a} deux disjointes
(par disjointes, on entend que les int\'{e}rieurs des connexions de
selles ne se rencontrent pas). Pour ce faire, on utilisera l'outil
g\'{e}om\'{e}trique suivant~:

\begin{Definition}
Soit $\Su$ une surface de translation. Un \emph{complexe} de
$\Su$ est une partie $K$ de $\Su$ qui est triangulable et dont
les bords sont des connexions de selles. On suppose de plus
que, si trois connexions de selles appartiennent \`{a} $K$ et
bordent un triangle ne contenant pas de singularit\'{e}, alors ce
triangle est inclus dans $K$.

Si $K$ est un complexe, on note $\Aire(K)$ sa surface, $|\partial K|$
le maximum des longueurs des connexions de selles formant le bord de
$K$, et $\Comp(K)$ la complexit\'{e} de $K$, i.e., le nombre de
connexions de selles n\'{e}cessaire pour trianguler $K$.
\end{Definition}

Notons que $\Comp(K)$ est bien d\'{e}fini puisque deux triangulations
diff\'{e}rentes de $K$ utilisent le m\^{e}me nombre de connexions de
selles, par la formule d'Euler. Comme toute triangulation de $K$ peut
\^{e}tre compl\'{e}t\'{e}e en une triangulation de $\Su$ (par
\ref{ExisteTriangulation}), et comme le nombre
d'ar\^{e}tes d'une triangulation de $\Su$ est constant, on en d\'{e}duit
que la complexit\'{e} d'un complexe est uniform\'{e}ment born\'{e}e par une
constante ne d\'{e}pendant que du genre de $\Su$ et du nombre de
singularit\'{e}s. \medskip

\noindent La proposition suivante permet de construire des complexes de
plus en plus grands (pour l'inclusion).

\begin{Proposition}
\label{prop:ConstruitComplexe}
Soient $K$ un complexe, et $\sigma$ une connexion de selles qui
traverse $\partial K$ ou qui ne rencontre pas $K$. Alors il existe un
complexe $K'$ contenant $K$ avec
  \begin{enumerate}
  \item $\Comp(K')> \Comp(K)$.
  \item $|\partial K'| \leq 2|\partial K| + |\sigma|$.
  \item $\Aire(K') \leq \Aire(K)+|\partial K|^2 + |\partial K|\cdot
  |\sigma|$.
  \end{enumerate}
\end{Proposition}
Notons que la proposition n'affirme \emph{pas} que $K'$ contient
$\sigma$ !
\begin{proof}
C'est de la g\'{e}om\'{e}trie \`{a} l'ancienne que n'aurait pas
d\'{e}savou\'{e}e Pythagore. On traite plusieurs cas. Dans tous les cas,
on construira $K'$ qui contient strictement $K$, la borne sur la
complexit\'{e} est donc triviale. Les autres bornes seront \'{e}galement
des cons\'{e}quences directes de la construction g\'{e}om\'{e}trique.

(1) Si $\sigma$ ne rencontre pas $K$ dans son int\'{e}rieur. On ajoute
simplement $\sigma$ \`{a} $K$ pour former $K'$ (et on remplit les
triangles \'{e}ventuels bord\'{e}s par $\sigma$ et deux connexions de
selles d\'{e}j\`{a} dans $K$).

(2) Si $\sigma$ rencontre $K$, mais une des extr\'{e}mit\'{e}s de $\sigma$
n'est pas dans $K$, appelons-la $P$. Partant de $P$, on parcourt
$\sigma$ jusqu'\`{a} rencontrer $K$, en un point $H$ qui appartient \`{a}
une connexion de selles $[AB]$ du bord de $K$, comme sur la figure
\ref{Augmente1}.

\begin{figure}[htbp]
\begin{center}
\psfrag{K}{$\scriptstyle K$} \psfrag{s}{$\scriptstyle \sigma$}
\psfrag{a}{$\scriptstyle A$} \psfrag{b}{$\scriptstyle B$}
\psfrag{h}{$\scriptstyle H$} \psfrag{p}{$\scriptstyle P$}

\includegraphics[width=11cm]{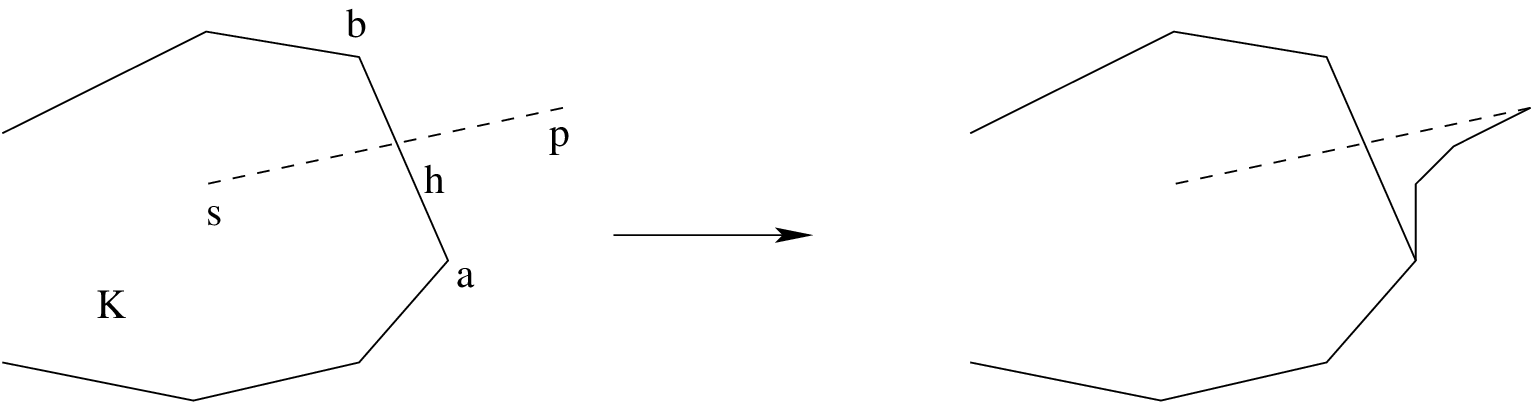}
\end{center}
\caption{\label{Augmente1}}
\end{figure}

On redresse le chemin $AHP$ petit \`{a} petit, jusqu'\`{a} rencontrer une
singularit\'{e}, puis on continue le redressement pour la portion
restante. On obtient une ligne bris\'{e}e reliant $A$ \`{a} $P$ qui
d\'{e}limite, avec $[AH]$ et $[HP]$, une zone de $\Su$ ne contenant pas
de singularit\'{e}. Si l'une des connexions de selles de la ligne bris\'{e}e
n'est pas dans $K$, on l'ajoute \`{a} $K$ pour former $K'$. Comme la
ligne bris\'{e}e est de longueur au plus $|AH|+|AP|$, $K'$ ainsi form\'{e}
v\'{e}rifie la conclusion de la proposition, i.e. 
$|\partial K'| \leq |\partial K| + |\sigma|$.

Si toutes les connexions de selles de la ligne bris\'{e}e sont dans $K$,
on r\'{e}p\`{e}te le processus de l'autre c\^{o}t\'{e}, entre $B$ et $P$. Le
cas qui pose probl\`{e}me est lorsque cette nouvelle ligne bris\'{e}e ne
contient que des connexions de selles dans $K$. Si l'une des lignes
bris\'{e}es entre $A$ et $P$, ou $B$ et $P$, n'est pas un segment, elle
contient un point interm\'{e}diaire $Q$. En reliant $Q$ \`{a} $A$ ou $B$,
on obtient une connexion de selles qui n'est pas dans $K$ (puisque $K$
ne rencontre pas $]HP[$), et en l'ajoutant \`{a} $K$ on obtient un
complexe $K'$ qui v\'{e}rifie les propri\'{e}t\'{e}s requises.

Il reste le cas o\`{u} $[AP]$ et $[BP]$ sont deux connexions de selles
appartenant \`{a} $K$. Mais alors, les connexions de selles $[AP]$,
$[BP]$ et $[AB]$ appartiennent \`{a} $K$ et bordent un triangle sans
singularit\'{e} dans son int\'{e}rieur. La condition que nous avons
impos\'{e}e dans la d\'{e}finition des complexes implique que $K$ contient
le triangle $ABP$, ce qui est une contradiction.

(3) Sinon, une partie de l'int\'{e}rieur de $\sigma$ n'est pas contenue
dans $K$, mais traverse de part et d'autre le bord de $K$, comme sur
la figure \ref{Augmente2}.

\begin{figure}[htbp]
\begin{center}
\psfrag{k}{$\scriptstyle K$} \psfrag{s}{$\scriptstyle \sigma$}
\psfrag{ou}{$\scriptstyle \textrm{ou}$}

\includegraphics[width=13cm]{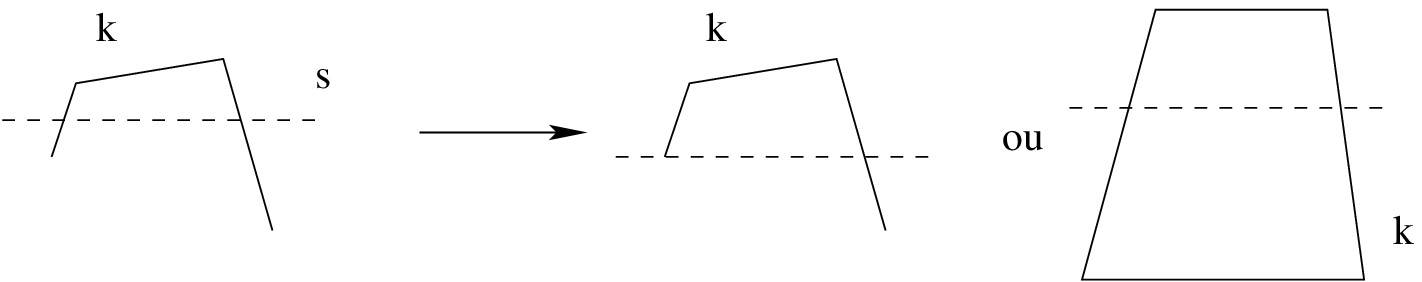}
\end{center}
\caption{\label{Augmente2}}
\end{figure}

On d\'{e}place alors $\sigma$ parall\`{e}lement \`{a} elle-m\^{e}me, dans une
direction, jusqu'\`{a} rencontrer une singularit\'{e}. Si on s'est
ramen\'{e} au cas (1) ou (2), on a gagn\'{e} (le nouveau $\sigma$,
translat\'{e}, est de longueur au plus $|\sigma|+|\partial K|$, et on
v\'{e}rifie que le $K'$ qu'on forme alors v\'{e}rifie bien les bornes requises).
Sinon, on a rencontr\'{e} une ligne
parall\`{e}le \`{a} $\sigma$, appartenant au bord de $K$ (contenant
\'{e}ventuellement des singularit\'{e}s dans son int\'{e}rieur, et
\'{e}ventuellement r\'{e}duite \`{a} un point). On essaie alors de
translater $\sigma$ dans l'autre direction, et soit on se ram\`{e}ne \`{a}
(1) ou (2), soit on a finalement une situation comme dans le dernier
cas de la figure \ref{Augmente2}. De plus, le quadrilat\`{e}re ainsi
form\'{e} n'est pas un triangle sans singularit\'{e} au bord (il
appartiendrait alors d\'{e}j\`{a} \`{a} $K$, par d\'{e}finition d'un
complexe). On va donc pouvoir relier une singularit\'{e} du bord
sup\'{e}rieur \`{a} une du bord inf\'{e}rieur, pour former $K'$, qui convient.
\end{proof}

\numero
On peut utiliser cette construction g\'{e}om\'{e}trique pour d\'{e}montrer
la proposition~\ref{prop:RencontrePasPetit}.

\begin{proof}[D\'{e}monstration de la proposition
\ref{prop:RencontrePasPetit}] Le c{\oe}ur de la preuve va \^{e}tre
de faire une r\'{e}currence sur la complexit\'{e} de certains
complexes. Comme la complexit\'{e} des complexes est uniform\'{e}ment
born\'{e}e, cette r\'{e}currence s'arr\^{e}tera \`{a} un certain moment. On
part de $\Su$ une surface de translation qui ne satisfait pas
la conclusion du th\'{e}or\`{e}me B.

Pour $i\in \N$, on note $(P_i)$ la propri\'{e}t\'{e} suivante : il existe
$\epsilon_n \to 0$, $T_n \to \infty$ et $S_n \subset
\Sbb^1$ avec $\inf_{n\in \N}\Leb( S_n) >0$ tels que, pour tout
$\theta \in S_n$, la surface $\Su_{\theta, T_n}$ admet un complexe
$K(n,\theta)$
de complexit\'{e} $\geq i$, avec $\Aire(K)\leq \epsilon_n$ et $|\partial
K|\leq \epsilon_n$.

On sait d\'{e}j\`{a} que $\Su$ satisfait la propri\'{e}t\'{e} $(P_1)$ : un
complexe de complexit\'{e} $1$ est simplement une connexion de selles,
et la proposition \ref{prop:RencontrePTpetit} donne donc la conclusion
de $(P_1)$.

Supposons maintenant que $(P_i)$ est satisfait, et consid\'{e}rons les
suites $\epsilon_n, T_n, S_n$ et $K(n,\theta)$ correspondantes. Soit
$c>0$. Si la
mesure de ``l'ensemble des $\theta$ tels que $\partial K(n,\theta)$ ne
rencontre pas de connexion de selles de longueur $\leq c$'' ne tend
pas vers $0$, on prend
une connexion de selles dans le bord de $K(n,\theta)$, elle v\'{e}rifie
la conclusion de la proposition \ref{prop:RencontrePasPetit} et on a
gagn\'{e}. Notons que $\partial K(n,\theta)$ est non vide pour $n$ assez
grand puisque l'aire de $K(n,\theta)$ tend vers $0$, ce qui garantit
que $K(n,\theta)$ ne recouvre pas toute la surface.

Sinon, d\'{e}montrons $(P_{i+1})$. Par un proc\'{e}d\'{e} diagonal (en
prenant une suite de valeurs de $c$ qui tend vers $0$), on
obtient une suite $\delta_n$ qui tend vers $0$, une sous-suite
$T'_n$ de $T_n$ (disons $T'_n=T_{j(n)}$) et des ensembles
$S'_n$ de mesure uniform\'{e}ment minor\'{e}e (inclus dans $S_{j(n)}$)
tels que, pour tout $\theta \in S'_n$, on ait un complexe
$L(n,\theta)=K(j(n),\theta)$ de complexit\'{e} au moins $i$ avec
$|\partial L|\leq \delta_n$ et $|\Aire(L)|\leq \delta_n$, ainsi
qu'une connexion de selles $\sigma(n,\theta)$ de longueur au
plus $\delta_n$ et qui rencontre le bord de $L$. En appliquant
la proposition \ref{prop:ConstruitComplexe} \`{a} $L$ et $\sigma$,
on obtient un complexe $K'$ de complexit\'{e} $\geq i+1$ d'aire et
de complexit\'{e} petites, born\'{e}es disons par
$\epsilon_n'=3\delta_n$. On a \'{e}tabli la propri\'{e}t\'{e} $(P_{i+1})$.

Si la surface $\Su$ ne v\'{e}rifiait pas la conclusion de la proposition
\ref{prop:RencontrePasPetit}, on d\'{e}montrerait ainsi $(P_1)$, puis
$(P_2)$, puis par r\'{e}currence $(P_i)$ pour tout $i$. Ceci est
impossible puisqu'il existe une borne uniforme sur la complexit\'{e} des
complexes de $\Su$. Ainsi, la proposition est d\'{e}montr\'{e}e.
\end{proof}

La contradiction viendra de la proposition suivante.
\begin{Proposition}
\label{prop:PetiteMesure}
Soient $\Su$ une surface de translation et $c>0$. Notons
$A(t,\epsilon,c)$ l'ensemble
des $\theta \in \Sbb^1$
tels que $\Su_{\theta,t}$ a une connexion de selles de longueur au
plus $\epsilon$ qui ne rencontre pas d'autre connexion de selles de
longueur au plus $c$. Il existe alors des
constantes
$\bar \epsilon>0, T>0$ et $K>0$ telles que, pour tout $t\geq T$, pour
tout $\epsilon \leq \bar\epsilon$,  $\Leb(A(t,\epsilon,c)) \leq K
\epsilon$.
\end{Proposition}
\begin{proof}
Posons $A=A(t,\epsilon,c)$, on veut estimer $\Leb(A)$. Soit $\theta
\in A$, on consid\`{e}re une connexion de selles $\alpha$ (dans $\Su$)
de longueur au plus $\epsilon$ et ne rencontrant pas de
connexion de selles de longueur au plus $c$ (dans
$\Su_{\theta,t}$). On note $I_\alpha$ l'ensemble des angles $\theta'$
tels que $\alpha$ est de longueur au plus $\epsilon$ dans
$\Su_{\theta',t}$, et $J_\alpha$ l'ensemble des angles $\theta'$ tels
que $\alpha$ est de longueur au plus $c$ dans
$\Su_{\theta',t}$. Ces ensembles sont des r\'{e}unions de deux
intervalles oppos\'{e}s sur le cercle $\Sbb^1$. Quitte \`{a} consid\'{e}rer
une seule des deux composantes (ou \`{a} projectiviser), on fera comme
si $I_\alpha$ et $J_\alpha$ \'{e}taient des intervalles. Soit $\hat
J_\alpha$ l'intervalle moiti\'{e} de $J_\alpha$ (i.e., l'intervalle de
m\^{e}me centre que $J_\alpha$ et de longueur moiti\'{e}).
Finalement, on note $K_\alpha$ l'ensemble des
angles $\theta'$ tels que $\alpha$ soit de longueur au plus $\epsilon$
dans $\Su_{\theta',t}$, et ne rencontre pas de connexion de selles de
longueur $\leq c$.
Si $\epsilon$ est assez petit et $t$ est assez grand, on a
  \begin{equation*}
  \{ \theta\} \subset K_\alpha \subset I_\alpha \subset \hat
  J_\alpha \subset J_\alpha.
  \end{equation*}
Un petit calcul trigonom\'{e}trique montre que, pour une certaine
constante $K$, $\Leb(I_\alpha) \leq K\epsilon \Leb(J_\alpha)$. On en
d\'{e}duit $\Leb(I_\alpha) \leq 2K \epsilon \Leb( \hat J_\alpha)$.

Notons $B$ l'ensemble des connexions de selles $\alpha$ telles que
$K_\alpha \not=\emptyset$. Ce sont les connexions de selles qui
servent \`{a} d\'{e}finir $A$. Autrement dit, $A=\bigcup_{\alpha \in B}
K_\alpha$.

Soient $\alpha,\beta \in B$ deux connexions de selles distinctes.
Si $\hat J_\alpha$ et $\hat J_\beta$
s'intersectent, alors $J_\alpha$ recouvre $\hat J_\beta$ ou
inversement. Dans le premier cas, comme $K_\beta \subset \hat
J_\beta$, on obtient un angle $\theta \in K_\beta \cap
J_\alpha$. Dans la surface de translation $\Su_{\theta,t}$, $\alpha$
est donc de longueur au plus $c$. Par d\'{e}finition de $K_\beta$, on en
d\'{e}duit que $\alpha$ et $\beta$ ne s'intersectent pas.

Plus g\'{e}n\'{e}ralement, si $\alpha_1,\dots, \alpha_k \in B$ sont deux
\`{a} deux distinctes et $\bigcap \hat J_{\alpha_i} \not=\emptyset$, on
montre de m\^{e}me que les connexions de selles $\alpha_i$ sont deux \`{a}
deux disjointes. Mais le nombre de connexions de selles deux \`{a} deux
disjointes est uniform\'{e}ment born\'{e} en fonction du genre de la
surface $\Su$ et du nombre de singularit\'{e}s (puisqu'on peut
compl\'{e}ter tout ensemble de connexions de selles disjointes pour
former une triangulation de la surface, qui a un nombre d'ar\^{e}tes
fix\'{e}). Il existe donc $p<\infty$ tel que chaque angle de $\Sbb^1$
appartienne \`{a} au plus $p$ intervalles $\hat J_\alpha$.

Finalement,
  \begin{align*}
  \Leb(A) & = \Leb( \bigcup_{\alpha\in B} K_\alpha)
  \leq \sum_{\alpha \in B} \Leb(I_\alpha)
  \leq 2 K \epsilon \sum_{\alpha\in B} \Leb( \hat J_\alpha)
  \\&
  \leq 2Kp \epsilon \Leb( \bigcup_{\alpha \in B} \hat J_\alpha)
  \leq 2Kp \epsilon \Leb(\Sbb^1).
  \qedhere
  \end{align*}
\end{proof}

\numero
On peut maintenant conclure :
\begin{proof}[D\'{e}monstration du th\'{e}or\`{e}me B]
Les propositions \ref{prop:RencontrePasPetit}
et \ref{prop:PetiteMesure} sont manifestement
incompatibles l'une avec l'autre, puisque l'ensemble $S_n$ donn\'{e} par
la premi\`{e}re proposition a une mesure born\'{e}e inf\'{e}rieurement,
alors que la seconde proposition montre que $\Leb(S_n)=O(\epsilon_n)
\to 0$.
\end{proof}

\section{Traduction de la preuve du th\'{e}or\`{e}me KMS en genre $1$}

\numero Dans ce paragraphe, on explicite diff\'{e}remment, avec un
point de vue plus g\'{e}om\'{e}trique, les objets qui sont apparus dans
la preuve du th\'{e}or\`{e}me KMS, dans le cas des tores. La notion
essentielle de ``d\'{e}part \`{a} l'infini'' utilis\'{e}e dans la preuve va
se reformuler en termes de ``l'espace des tores plats''. Notons
que, pour pouvoir parler de systole, il faut avoir au moins une
singularit\'{e} (\'{e}ventuellement artificielle, i.e., avec $k_i=0$).
On travaillera donc sur des tores ayant un point marqu\'{e} (plus
ou moins arbitraire).

\numero Soit $(e_1,e_2)$ une base directe de $\R^2$. On peut
alors former une surface de translation de genre $1$, avec un
point marqu\'{e}, comme suit : on part du polygone $P=\{x e_1 +y
e_2,\ 0 \leq x \leq 1, 0 \leq y \leq 1\}$ puis on recolle ses
c\^{o}t\'{e}s oppos\'{e}s parall\`{e}les, comme en~\ref{par:RecollePolygones}.
On obtient ainsi une surface de translation $\Su$, sur laquelle
on marque l'image de l'origine. De plus, tout tore plat avec un
point marqu\'{e} peut s'obtenir de cette mani\`{e}re.

Soit $\T_1$ l'ensemble des bases directes de $\R^2$, cet
ensemble fournit donc une param\'{e}trisation de ``l'espace des
tores plats avec un point marqu\'{e}''.

Notons $e_1=\left(\begin{smallmatrix}a \\ b
\end{smallmatrix}\right)$ et $e_2=\left(
\begin{smallmatrix}c \\ d \end{smallmatrix}\right)$.
On a une identification entre $\T_1$ et
$\textrm{GL}^{+}_2(\R)$ donn\'{e}e par $(e_1,e_2) \mapsto
A_{(e_1,e_2)}= \left(\begin{smallmatrix}a & b \\ c & d
\end{smallmatrix}\right)$ (le fait d'avoir transpos\'{e} dans cette
identification nous servira plus loin). Si $\T^{(1)}_1$ d\'{e}signe
l'ensemble des tores d'aire $1$ alors $\T^{(1)}_1 = \Sl$ via
l'identification pr\'{e}c\'{e}dente. \medskip

\noindent Analysons l'action de $\Sl$, d\'{e}crite dans la
section~\ref{sec:action}, sur $\T_1$. Si $M\in \Sl$ alors
l'image par $M$ du tore plat donn\'{e} par $(e_1,e_2)$ est le tore
plat donn\'{e} par $(Me_1,Me_2)$. Via l'identification ci-dessus,
$M$ agit donc sur $\textrm{GL}^{+}_2(\R)$ comme la
multiplication \`{a} droite par la matrice transpos\'{e}e $M^{t}$.

En particulier, l'action du flot g\'{e}od\'{e}sique sur $\T^{(1)}_1$
correspond \`{a} l'action par multiplication \`{a} droite de $g_t$ sur
$\Sl$.

\numero L'ensemble $\T_1$ d\'{e}crit ci-dessus n'est pas tr\`{e}s
satisfaisant. En effet, deux bases directes $(e_1,e_2)$ et
$(e'_1,e'_2)$ peuvent engendrer deux surfaces de translation
$\Su$ et $\Su'$ isomorphes (i.e., il existe un diff\'{e}omorphisme
de $\Su$ dans $\Su'$ dont la diff\'{e}rentielle dans l'atlas de
translation est partout \'{e}gale \`{a} l'identit\'{e}). Notons $\mathcal
M_1$ le quotient de $\T_1$ par cette relation d'\'{e}quivalence.

On va d\'{e}crire $\mathcal M_1$ via l'identification ci-dessus.
Deux bases directes $(e_1,e_2)$ et $(e'_1,e'_2)$ sont
\'{e}quivalentes si et
seulement si il existe une matrice $B=\left(\begin{smallmatrix}x & y\\
z & t
\end{smallmatrix}\right)$ dans $\textrm{SL}_2(\Z)$
telle que $e'_1=xe_1+ye_2$ et $e'_2=ze_1+te_2$, i.e., si
$(e_1,e_2)$ et $(e'_1,e'_2)$ sont deux bases d'un m\^{e}me r\'{e}seau.
Cette condition se lit encore $A_{(e'_1,e'_2)}=B
A_{(e_1,e_2)}$. Ainsi, $\mathcal M_1$ est identifi\'{e} au quotient
(\`{a} gauche) $\textrm{SL}_2(\Z) \setminus \textrm{GL}^{+}_2(\R)$.
Soit $\mathcal M_1^{(1)}$ le sous-ensemble de $\mathcal M_1$
form\'{e} des surfaces d'aire $1$, il est identifi\'{e} \`{a}
$\textrm{SL}_2(\Z)\setminus \Sl$.

\noindent Notons que l'action de $g_t$ passe au quotient sur $\mathcal
M^{(1)}_1$ (une action \`{a} gauche et une action \`{a} droite commutent
toujours) et s'interpr\`{e}te alg\'{e}briquement comme la
multiplication \`{a} droite de $g_t$ sur $\textrm{SL}_2(\Z) \setminus
\Sl$.

\numero 
\label{prop:recurrent}
On peut alors relier la notion na\"{\i}ve de d\'{e}part \`{a}
l'infini donn\'{e}e en \ref{def:PartirInfini} avec la topologie de
$\textrm{SL}_2(\Z) \setminus \Sl$.

\begin{NoNumberProposition}
\label{DepartInfiniGenre1} Un tore $\mathbb T^2 \in
\textrm{SL}_2(\Z) \setminus \Sl$ part \`{a} l'infini (au sens de la
d\'{e}finition~\ref{def:PartirInfini}) si et seulement si $g_t
\mathbb T^2$ quitte tout compact de $\textrm{SL}_2(\Z)
\setminus \Sl$ quand $t\to \infty$.
\end{NoNumberProposition}

\begin{proof}
Une suite de tores dans $\textrm{SL}_2(\Z) \setminus \Sl$
quitte tout compact si et seulement si la systole tend vers $0$
le long de cette suite. La proposition en d\'{e}coule
imm\'{e}diatement.
\end{proof}

Ainsi, le th\'{e}or\`{e}me A se comprend bien en termes d'espace des
r\'{e}seaux.

\numero On aurait aussi pu d\'{e}crire $\T_1$ comme l'ensemble des
atlas plats sur un tore, quotient\'{e} par l'ensemble des
diff\'{e}omorphismes qui sont isotopes \`{a} l'identit\'{e}. De ce point de
vue g\'{e}om\'{e}trique $\mathcal M^{(1)}_1$ devient alors l'ensemble
des atlas plats d'aire $1$ sur un tore, quotient\'{e} par
l'ensemble des diff\'{e}omorphismes. Cela revient aussi \`{a} prendre
$\T^{(1)}_1$ modulo le groupe modulaire, i.e., le groupe de
tous les diff\'{e}omorphismes (pr\'{e}servant l'orientation) modulo
ceux isotopes \`{a} l'identit\'{e}. En effet si $\phi$ est un
diff\'{e}omorphisme du tore $\mathbb T^2$ et $\{(U_i,z_i)\}_i$ est
un atlas plat, alors $\{(\phi(U_i),z_i\circ \phi^{-1})\}_i$ est
un autre atlas plat repr\'{e}sentant le m\^{e}me tore plat. Ici le
groupe modulaire n'est autre que $\textrm{SL}_2(\Z)$. Ce point
de vue va \^{e}tre utile en genre sup\'{e}rieur.

\numero Nous allons maintenant red\'{e}montrer le th\'{e}or\`{e}me B, en
genre $1$. Pour cela, il est utile d'introduire un autre point
de vue sur l'espace $\textrm{SL}_2(\Z) \setminus \Sl$, celui de
la g\'{e}om\'{e}trie hyperbolique.

Le groupe $\Sl$ agit sur le demi-plan de Poincar\'{e} 
$\mathbb H=\{z\in \C,\ \textrm{Im}(z) > 0 \}$, par

  \begin{equation}
  \left(\begin{matrix}
  a & b\\c&d
  \end{matrix}\right)\cdot z= \frac{az+b}{cz+d}.
  \end{equation}
Cette action pr\'{e}serve la m\'{e}trique hyperbolique. De plus, cette 
action est transitive et le 
stabilisateur du point $i\in \mathbb H$ est $\textrm{SO}(2)$,
donc $\mathbb H \simeq \Sl / \textrm{SO}(2)$, et $\Sl$ est
identifi\'{e} au fibr\'{e} unitaire tangent \`{a} $\mathbb H$.

\begin{figure}[htbp]
\begin{center}
\psfrag{H}{$\scriptstyle \mathbb H$}  \psfrag{+}{$\scriptstyle
\frac{1}{2}$} \psfrag{-}{$\scriptstyle -\frac{1}{2}$}
\psfrag{0}{$\scriptstyle 0$} \psfrag{i}{$\scriptstyle i$}
\psfrag{j}{$\scriptstyle j$} \psfrag{infty}{$\scriptstyle \infty$}
\psfrag{-1}{$\scriptstyle -1$} \psfrag{1}{$\scriptstyle 1$}

\includegraphics[width=8cm]{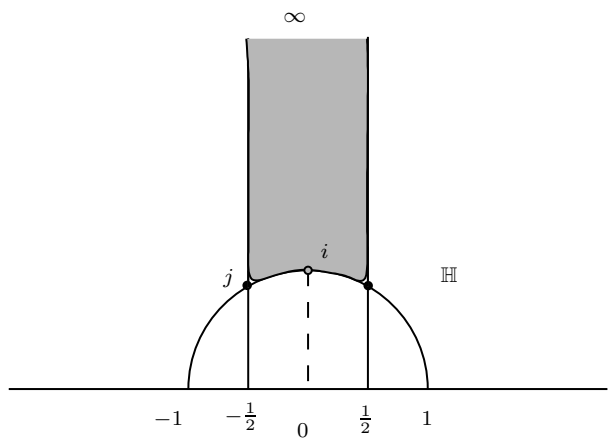}
\end{center}
\caption{ \label{fig:surface:modulaire} Un domaine fondamental pour
l'action de $\textrm{SL}_2(\Z)$ sur $\mathbb H$.
}
\end{figure}

Un domaine fondamental pour l'action (\`{a} gauche) de
$\textrm{SL}_2(\Z)$ sur $\mathbb H$ est pr\'{e}sent\'{e}
figure~\ref{fig:surface:modulaire}. La surface modulaire est
par d\'{e}finition le quotient $\textrm{SL}_2(\Z)\setminus \mathbb H$.
Cette surface est isomorphe \`{a} une sph\`{e}re priv\'{e}e d'un point avec
deux singularit\'{e}s pour la m\'{e}trique hyperbolique. Notons que les
deux singularit\'{e}s correspondent au tore plat carr\'{e} et au tore
plat hexagonal.

L'espace des surfaces de translation de genre $1$ avec un point
marqu\'{e}, i.e. $\textrm{SL}_2(\Z) \setminus \Sl$, est donc
identifi\'{e} au fibr\'{e} unitaire tangent \`{a} la surface modulaire. De
plus, l'action de $g_t$ correspond au flot g\'{e}od\'{e}sique sur la
surface modulaire.

\numero \label{surface:veech} On peut maintenant red\'{e}montrer en
genre $1$ une version forte du th\'{e}or\`{e}me B, en utilisant des
propri\'{e}t\'{e}s bien connues du flot g\'{e}od\'{e}sique sur la surface
modulaire.

\begin{NoNumberProposition}
Soit $\Su$ une surface de translation de genre $1$ avec un
point marqu\'{e}. Alors
  \begin{equation}
  \label{eq:PartInfiniGenre1}
  \{\theta\in \Sbb^1,\ \Su \textrm{ part \`{a} l'infini
  dans la direction }\theta\}
  \end{equation}
est au plus d\'{e}nombrable.
\end{NoNumberProposition}
\begin{proof}
On utilisera la propri\'{e}t\'{e} suivante du flot g\'{e}od\'{e}sique sur la
surface modulaire : il existe un compact $K$  tel que toute
g\'{e}od\'{e}sique qui ne part pas verticalement dans le cusp revient
une infinit\'{e} de fois dans $K$ (i.e. les g\'eod\'esiques qui ne 
partent pas vers le cusp sont r\'ecurrentes, c'est la proposition 
du point~\ref{prop:recurrent}).

Soit $\Su$ une surface de translation avec un point marqu\'{e},
consid\'{e}rons $x\in \textrm{SL}_2(\Z)\setminus \Sl$ le point
correspondant, et $y=\pi(x)$ son image dans la surface
modulaire $\textrm{SL}_2(\Z)\setminus \mathbb H$. L'ensemble
\eqref{eq:PartInfiniGenre1} est alors en bijection avec les
g\'{e}od\'{e}siques issues de $y$ qui partent \`{a} l'infini dans la
surface modulaire.

Pour tout $n>0$, l'ensemble $G_n$ des g\'{e}od\'{e}siques partant de
$y$ et qui sont verticales dans le cusp apr\`{e}s le temps $n$ est
discret, donc fini par compacit\'{e} de $\Sbb^1$. L'ensemble
\eqref{eq:PartInfiniGenre1}, en bijection avec $\bigcup G_n$,
est donc d\'{e}nombrable, ce qui conclut la preuve.
\end{proof}

Cette preuve ne s'adapte pas en genre sup\'{e}rieur (et le r\'{e}sultat
de d\'{e}nombrabilit\'{e} devient faux !) car les propri\'{e}t\'{e}s
hyperboliques du flot $g_t$ deviennent beaucoup plus faibles
qu'en genre $1$. La morale du th\'{e}or\`{e}me B est qu'il reste
cependant suffisamment d'hyperbolicit\'{e} pour obtenir des
r\'{e}sultats pr\'{e}cis.

\section{L'espace des modules des surfaces de translation}
\label{sec:espace:module}

\numero Comme nous venons de le voir dans la section
pr\'{e}c\'{e}dente, le bon cadre conceptuel pour d\'{e}montrer le th\'{e}or\`{e}me
KMS est d'introduire une structure de ``vari\'{e}t\'{e}'' (plus
pr\'{e}cis\'{e}ment d'orbifold) sur ``l'espace des surfaces de
translation'' (en genre fix\'{e} et avec une combinatoire \'{e}galement
fix\'{e}e). Pour que cette notion soit utile, on voudrait avoir un
analogue de la proposition~\ref{DepartInfiniGenre1}, i.e., le
d\'{e}part \`{a} l'infini au sens de la
d\'{e}finition~\ref{def:PartirInfini} doit correspondre au fait de
quitter tout compact de l'espace des surfaces de translation.
\medskip

\noindent Dans cette section, sans rien d\'{e}montrer, on va d\'{e}crire
une telle construction.

\numero Fixons $\Su$ une surface compacte connexe, de genre
$g$. Fixons aussi $\Sigma = \{P_1,\dots,P_n\}$ un sous-ensemble
fini de $\Su$, et $\kappa=(k_1,\dots, k_n)$ des entiers tels
que $\sum (k_i+1)=2g-2$. \medskip

\noindent On note $\A(\Su, \kappa)$ l'ensemble des structures
de translation sur $X$ ayant des singularit\'{e}s coniques d'angles
$2\pi(k_i+1)$ en $P_i$. C'est un ensemble d'atlas, il est
\'{e}norme (et en particulier beaucoup trop gros pour qu'on puisse
travailler avec lui). On peut n\'{e}anmoins le munir naturellement
d'une topologie : si $\omega \in \A(\Su, \kappa)$ et
$\epsilon>0$, on d\'{e}finit le $\epsilon$-voisinage de $\omega$
comme \'{e}tant l'ensemble des atlas de translation $\omega'$ sur
$\Su$ tels que, pour tout $x\in \Su \moins \Sigma$, les
compos\'{e}es $z\circ {z'}^{-1}$ et $z'\circ z^{-1}$ de cartes de
translation autour de $x$ (correspondant \`{a} $\omega$ et
$\omega'$) ont des d\'{e}riv\'{e}es $\epsilon$-proches de l'identit\'{e}.
\medskip

\noindent Intuitivement, $\omega'$ est dans le
$\epsilon$-voisinage de $\omega$ si $\omega'$ et $\omega$ sont
d'accord \`{a} $\epsilon$ pr\`{e}s sur la direction des vecteurs dans
$\R^2$, leur norme, les horizontales, les verticales, etc.

\numero
Soit $\T_g(\Su,\kappa)$ obtenu en identifiant deux \'{e}l\'{e}ments $\omega$ et
$\omega'$ de $\A(\Su,\kappa)$ s'il existe un hom\'{e}omorphisme de $\Su$
fixant les $P_i$, isotope \`{a} l'identit\'{e} (relativement aux $P_i$) et
envoyant $\omega$ sur $\omega'$. On munit $\T_g(\Su,\kappa)$ de la topologie
quotient. C'est ``l'espace de Teichm\"{u}ller'' des surfaces de translation. \medskip

\noindent En g\'{e}n\'{e}ral, il faut se m\'{e}fier de la topologie quotient, mais ici
tout va bien : l'espace $\T_g(\kappa)$ est s\'{e}par\'{e}, m\'{e}trisable, localement
compact, et m\^{e}me localement hom\'{e}omorphe \`{a} un certain espace
$\C^d$ (ici $d$ est $2g+n-1$ comme nous allons le voir).

\numero On peut m\^{e}me aller plus loin et munir $\T_g(\Su, \kappa)$
d'une structure de vari\'{e}t\'{e} (analytique), comme suit. Si
$\omega\in \A(\Su,\kappa)$, on peut utiliser la structure de
translation pour relever les chemins sur $\Su$, en partant de
$0$. L'extr\'{e}mit\'{e} du rel\`{e}vement est invariante par homotopie, et
on obtient donc une application de $H_1(\Su, \Sigma ; \Z)$ dans
$\R^2$, i.e., un \'{e}l\'{e}ment du groupe de cohomologie
$H^1(\Su,\Sigma; \R^2)$. Cet \'{e}l\'{e}ment ne change pas si l'on
remplace $\omega$ par une structure de translation qui lui est
isotope. En passant au quotient, on obtient donc une
application canonique
  \begin{equation*}
  \Theta : \T_g(\Su,\kappa) \to H^1(\Su, \Sigma ; \R^2).
  \end{equation*}

\begin{Proposition}
Pour la topologie d\'{e}finie ci-dessus sur $\T_g(\Su, \kappa)$,
l'application $\Theta$ est un hom\'{e}omorphisme local.
\end{Proposition}

\noindent On peut donc utiliser $\Theta$ comme carte locale, pour mettre une
structure de vari\'{e}t\'{e} sur $\T_g(\kappa)$.
Dans le cas du genre $1$, $\T_1$ correspond simplement \`{a}
$\textrm{GL}^{+}_2(\R)$ (avec sa structure de vari\'{e}t\'{e} usuelle). 
Si $g \geq 2$, on obtient bien que $\T_g(\kappa)$ est une varit\'et\'e de dimension 
(r\'eelle) $2(2g+n-1)$.

\numero
Pour obtenir un espace plus petit, avec de meilleures propri\'{e}t\'{e}s
de compacit\'{e}, et donc de r\'{e}currence pour la dynamique,
il faut quotienter encore plus. Soit $\M_g(\kappa)$ le quotient
de $\T_g(\kappa)$ par l'action du groupe modulaire. Autrement dit, on identifie
deux \'{e}l\'{e}ments $\omega$ et $\omega'$ de $\A$ si $\omega'$ s'obtient \`{a} partir de
$\omega$ par un diff\'{e}omorphisme de $\Su$ fixant les $P_i$ (mais pas
n\'{e}cessairement isotope \`{a} l'identit\'{e}). \medskip

\noindent L'action du groupe modulaire sur $\T_g(\kappa)$ n'est pas
libre, mais presque. Ainsi, $\M_g(\kappa)$ n'est pas naturellement muni
d'une structure de vari\'{e}t\'{e} (il a des singularit\'{e}s), mais
presque : c'est un \emph{orbifold}, i.e., il est localement
diff\'{e}omorphe \`{a} un espace $\C^d$ sauf en un nombre fini de
points. On notera aussi $\M_g^{(1)}(\kappa)$ l'ensemble des surfaces dans
$\M_g(\kappa)$ d'aire $1$.

\noindent Dans ce cadre, le r\'{e}sultat de compacit\'{e} na\"{\i}f qu'on a
d\'{e}montr\'{e} dans le corollaire \ref{cor:Compacite} se traduit de la
mani\`{e}re suivante. Notons $\sys(\omega)$ la systole de $\omega\in \M_g^{(1)}(\kappa)$.

\begin{Theorem}
Pour tout $\epsilon>0$, l'ensemble des $\omega \in \M_g^{(1)}(\kappa)$ tels que
$\sys(\omega) \geq \epsilon$ est compact.
\end{Theorem}

\noindent Ainsi, une suite $\omega_n$ sort de tout compact si et seulement si
$\sys(\omega_n) \to 0$. La d\'{e}finition de d\'{e}part \`{a} l'infini donn\'{e}e
en \ref{def:Convergence} co\"{\i}ncide donc avec celle d\'{e}coulant de la
topologie naturelle sur l'espace des modules.

\numero
L'espace $H^1(\Su, \Sigma ; \R^2)$ est muni d'une mesure de Lebesgue
canonique (donnant covolume $1$ au r\'{e}seau $H^1(\Su, \Sigma ;
\Z^2)$). En la tirant en arri\`{e}re par $\Theta$, on obtient une mesure
canonique sur $\T_g(\kappa)$. Cette mesure est invariante sous l'action du
groupe modulaire, et passe donc au quotient sur $\M_g(\kappa)$. Elle induit
m\^{e}me une mesure sur $\M_g^{(1)}(\kappa)$.

\numero Un r\'{e}sultat important de Masur et Veech
(voir~\cite{Masur:82,Veech:82}) affirme que cette mesure est de
masse finie sur $\M_g^{(1)}(\kappa)$. De plus elle pr\'{e}serv\'{e}e par le
flot $g_t$. Masur et Veech montrent m\^eme que le flot $g_t$ est ergodique sur 
toute composante connexe de $\M_g^{(1)}(\kappa)$ (qui a au plus trois 
composantes connexes). \\
Quand on parle de ``presque toute
surface de translation'', on fait r\'{e}f\'{e}rence \`{a} cette mesure
canonique.

\numero \label{surface:reseau} Un cas particuli\`{e}rement
int\'{e}ressant est fourni lorsque le stabilisateur
  \begin{equation*}
  \textrm{SL}(\Su,\omega):=\{A\in \Sl,\ A\cdot (\Su,\omega) \simeq (\Su,\omega)\}
  \end{equation*}
d'une surface de translation $(\Su,\omega)$ est un r\'{e}seau de
$\textrm{SL}_2(\R)$. On dit alors que $(\Su,\omega)$ est une
{\it surface de Veech}. L'id\'{e}e de la preuve donn\'{e}e
en~\ref{surface:veech} peut alors s'adapter (la propri\'{e}t\'{e}
essentielle que nous avons utilis\'{e}e est le fait que
$\textrm{SL}_2(\Z)$ est un r\'{e}seau dans $\textrm{SL}_2(\R)$). On
peut alors obtenir des propri\'{e}t\'{e}s remarquables sur le flot
directionnel sur $\Su$ (voir~\cite{Veech:89}). Notamment, pour
une surface de Veech, les directions minimales pour le flot
directionnel sont {\it exactement} les directions uniquement
ergodiques.

\section{Raffinements, r\'{e}sultats suppl\'{e}mentaires}
\label{sec:raffinements}

Dans ce paragraphe, on discute les raffinements possibles aux
th\'{e}or\`{e}mes~$A$ et~$B$.

\numero
Tout d'abord, pour le th\'{e}or\`{e}me A, on a \'{e}quivalence entre non-unique
ergodicit\'{e} et d\'{e}part \`{a} l'infini, en genre $1$. Cela reste vrai pour
les surfaces de Veech (voir~\ref{surface:reseau}) (\cite{Veech:89}).
Mais ce n'est plus vrai en g\'{e}n\'{e}ral. Mentionnons les r\'{e}sultats suivants.
\begin{enumerate}
\item Pour toute fonction $f(t)$ qui tend vers $0$ quand $t\to \infty$,
il existe une surface $\Su$ telle que $\sys(g_t \Su)$ tend vers $0$
quand $t \to \infty$, mais $f(t)=o ( \sys(g_t \Su))$ (voir~\cite{Cheung:slow}).
Autrement dit, on peut converger arbitrairement lentement vers l'infini.
\item Pour toute composante connexe de l'espace des modules, il existe
une constante $c>0$ telle que, si $t^{-c}=o( \sys(g_t \Su))$ alors le
flot vertical sur $\Su$ est uniquement ergodique. Autrement dit, si on
tend vers l'infini mais suffisamment lentement, alors la conclusion du
th\'{e}or\`{e}me~$A$ reste valable. Ce th\'{e}or\`{e}me n'est pas vide puisqu'il
existe effectivement des g\'{e}od\'{e}siques qui partent lentement vers
l'infini d'apr\`{e}s le point pr\'{e}c\'{e}dent.
\item Pour toute composante connexe de l'espace des modules, il existe
une constante $c' >c$ et une surface $\Su$ telle que
$t^{-c'}=o(\sys(g_t\Su))$, mais le flot vertical sur $\Su$ n'est pas
uniquement ergodique. Ainsi, la vitesse de d\'{e}croissance en $t^{-c}$
est vraiment la ``vitesse critique''.

Ces deux derniers points sont dus \`{a} Cheung et Eskin
(voir~\cite{Cheung:Eskin}). On pourra aussi consulter~\cite{Cheung:Masur}.

\item Il existe une version combinatoire du th\'{e}or\`{e}me A
(voir~\cite{Boshernitzan}). \\
Une section de Poincar\'{e} du flot $\mathcal F_\theta$ sur un
intervalle fournit un \'{e}change d'intervalles $T$. Les points de
discontinuit\'{e}s correspondent aux s\'{e}paratrices. Notons $m_n$ la
longueur du plus petit intervalle \'{e}chang\'{e} par $T^n$. Alors, si
$T$ n'est pas uniquement ergodique, $nm_n \rightarrow \infty$
quand $n\rightarrow \infty$.

\end{enumerate}

\numero Pour le th\'{e}or\`{e}me~$B$, on peut se demander quelle est
la taille pr\'{e}cise de l'ensemble
  \begin{equation*}
   \{\theta\in \Sbb^1,\ \Su \textrm{ part \`{a} l'infini
  dans la direction }\theta\}.
  \end{equation*}
On sait d\'{e}j\`{a} qu'il est de mesure nulle. En fait on peut dire un peu plus sur
le sous-ensemble suivant~:
\begin{multline}
\label{equation:veech}
\{\theta\in \Sbb^1,\textrm{ le flot } \mathcal F_\theta
\textrm{ sur } \Su \textrm{ n'est pas uniquement ergodique} \} \subset \\
\{\theta\in \Sbb^1,\ \Su \textrm{ part \`{a} l'infini dans la direction }\theta\}.
\end{multline}
On peut chercher \`{a} estimer sa dimension de Hausdorff. On sait
d\'{e}j\`{a} qu'en genre $1$, l'ensemble des directions exceptionnelles
co\"{\i}ncide avec les rationnels et est donc d\'{e}nombrable. Ce n'est
plus vrai en genre plus grand, mais~:
\begin{enumerate}
\item Pour toute surface de translation, la dimension de Hausdorff de
cet ensemble est au plus $1/2$ (voir~\cite{Masur:HD}).
\item Il existe des surfaces de translation pour lesquelles cette
dimension est exactement $1/2$~\cite{Cheung} (ces exemples proviennent
de billards carr\'{e}s avec un mur, voir figure~\ref{fig:exemple:4:slit}).
\item Pour toute composante connexe de l'espace des modules $\M_g(\kappa)$, il
existe une constante $\delta \in ]0,1/2]$ telle que, pour presque
toute surface de translation dans cette composante connexe, la
dimension en question est \'{e}gale \`{a} $\delta$~\cite{Masur:Smillie:1}.
\item Pour les surfaces de Veech (voir~\ref{surface:reseau})
l'inclusion~(\ref{equation:veech}) est une \'{e}galit\'{e} (\cite{Veech:89}).
\end{enumerate}


\bibliography{biblio}

\begin{thebibliography}{KMS86}


\bibitem[Arn88]{Arnoux}
{\scshape P.~Arnoux} -- {\og
Ergodicit\'e g\'en\'erique des billards polygonaux 
(d'apr\`es Kerckhoff, Masur, Smillie)\fg}, S\'eminaire Bourbaki
\emph{Ast\'erisque}, \textbf{696} Vol. 1987/88.  No.~161-162  (1988), 
p.~203--221.


\bibitem[Bos85]{Boshernitzan}
{\scshape M.~BoshernitzanY} -- {\og
A condition for minimal interval exchange maps to be uniquely ergodic\fg}, 
\emph{Duke Math. J.}, \textbf{52} (1985), no.~3, p.~723--752.


\bibitem[CE06]{Cheung:Eskin}
{\scshape Y.~Cheung {\normalfont \smfandname} A.~Eskin} -- {\og
Unique Ergodicity of Translation Flows \fg}, 
\emph{Proceedings of the Partially hyperbolic dynamics, laminations, 
and Teichm\"uller flow Workshop}, Toronto, Jan. 5-9, (2006), \`a para\^itre.


\bibitem[CM06]{Cheung:Masur}
{\scshape Y.~Cheung {\normalfont \smfandname} H.~Masur} -- {\og A divergent 
Teichm\"uller geodesic with uniquely ergodic vertical foliation\fg}, 
\emph{Israel J. Math.} \textbf{152} (2006), p.~1--15.

\bibitem[Che03]{Cheung}
{\scshape Y.~Cheung} -- {\og Hausdorff dimension of the set of nonergodic
  directions\fg}, \emph{Ann. of Math. (2)} \textbf{158} (2003), no.~2,
  p.~661--678, With an appendix by M. Boshernitzan.

\bibitem[Che04]{Cheung:slow}
\bysame , {\og Slowly divergent geodesics in moduli space\fg}, \emph{Conform.
  Geom. Dyn.} \textbf{8} (2004), p.~167--189 (electronic).


\bibitem[FK36]{FoKe}
{\scshape R.~H. Fox {\normalfont \smfandname} R.~B. Kershner} -- {\og
  Concerning the transitive properties of geodesics on a rational
  polyhedron\fg}, \emph{Duke Math. J.} \textbf{2} (1936), no.~1, p.~147--150.

\bibitem[HM79]{Hubbard:Masur}
{\scshape J.~Hubbard {\normalfont \smfandname} H.~Masur} -- {\og Quadratic
  differentials and foliations\fg}, \emph{Acta Math.} \textbf{142} (1979),
  no.~3-4, p.~221--274.

\bibitem[Kea75]{Keane}
{\scshape M.~Keane} -- {\og Interval exchange transformations\fg}, 
\emph{Math. Z.} \textbf{141} (1975), no.~3-4, p.~25--31.


\bibitem[KMS86]{Kerckhoff:Masur:Smillie}
{\scshape S.~Kerckhoff, H.~Masur {\normalfont \smfandname} J.~Smillie} -- {\og
  Ergodicity of billiard flows and quadratic differentials\fg}, \emph{Ann. of
  Math. (2)} \textbf{124} (1986), no.~2, p.~293--311.


\bibitem[Mas82]{Masur:82}
{\scshape H.~Masur} -- {\og Interval exchange transformations and measured foliations\fg}, 
\emph{Ann. of Math. (2)} \textbf{115} (1982), no.~1, p.~169--200.


\bibitem[Mas92]{Masur:HD}
\bysame , {\og Hausdorff dimension of the set of nonergodic
  foliations of a quadratic differential\fg}, \emph{Duke Math. J.} \textbf{66}
  (1992), no.~3, p.~387--442.


\bibitem[Mas06]{Masur:06}
\bysame , {\og Ergodic theory of translation surfaces\fg}, 
\emph{Handbook of dynamical systems} 
Vol. 1B, (2006), Amsterdam Elsevier B.~V., p.~527--547.


\bibitem[MT02]{Masur:Tabachnikov}
{\scshape H.~Masur {\normalfont \smfandname} S.~Tabachnikov} -- {\og 
Rational billiards and flat structures\fg}, \emph{Handbook of dynamical systems} 
Vol. 1A, (2002), North-Holland, Amsterdam, p.~1015--1089.


\bibitem[MS91]{Masur:Smillie:1}
{\scshape H.~Masur {\normalfont \smfandname} J.~Smillie} -- {\og Hausdorff
  dimension of sets of nonergodic measured foliations\fg}, \emph{Ann. of Math.
  (2)} \textbf{134} (1991), no.~3, p.~455--543.

\bibitem[Vee68]{Veech:68}
{\scshape W.~Veech} -- {\og The equicontinuous structure relation for minimal 
Abelian transformation groups\fg}, \emph{Amer. J. Math.} 
\textbf{90} (1968), p.~723--732.


\bibitem[Vee82]{Veech:82}
\bysame , {\og Gauss measures for transformations on the 
space of interval exchange maps\fg}, \emph{Ann. of Math. (2)} 
\textbf{115} (1982), no.~1, p.~201--242.


\bibitem[Vee89]{Veech:89}
\bysame , {\og Teichm\"uller curves in moduli space, Eisenstein series and an application to 
triangular billiards\fg}, \emph{Invent. Math.} 
\textbf{97} (1989), no.~3, p.~553--583.


\bibitem[ZK75]{KaZa}
{\scshape A.~Zemljakov {\normalfont \smfandname} A.~Katok} -- {\og Topological
  transitivity of billiards in polygons\fg}, \emph{Mat. Zametki} \textbf{18}
  (1975), no.~2, p.~291--300.


\bibitem[Zor06]{Zorich:06}
{\scshape A.~Zorich} -- {\og Flat surfaces\fg}, 
\emph{collection "Frontiers in Number Theory, Physics and Geometry, 
Volume 1: On random matrices, zeta functions and dynamical systems},
P.~Cartier, B.~Julia, P.~Moussa, P.~Vanhove (Editors),
Springer-Verlag, Berlin, (2006), p.~439-586.

\end{thebibliography}

\end{document}